\documentclass[lettersize,journal]{IEEEtran}
\usepackage{amsmath,amsthm,amsfonts,amssymb,amscd, amsxtra,color}
\usepackage[active]{srcltx}
\usepackage{algorithm}
\usepackage{algpseudocode}
\usepackage{url}
\usepackage{graphicx}
\usepackage{cite}
\usepackage{lscape}
\usepackage{subfig}
\usepackage{multirow}
\theoremstyle{definition}
\newtheorem{theorem}{Theorem}
\newtheorem{lemma}{Lemma}
\newtheorem{definition}{Definition}
\newtheorem{corollary}{Corollary}
\newtheorem{proposition}{Proposition}

\newcommand{\beq}{\begin{equation}}
\newcommand{\eeq}{\end{equation}}
\newcommand{\tp}{^\top}

\newcommand{\lf}{\left}
\newcommand{\rg}{\right}

\def\min{\operatorname{min}}
\def\max{\operatorname{max}}

\newcommand{\R}{\mathbb R}

\graphicspath{{figs/}}

\begin{document}

\title{A Relaxation Method for Nonsmooth Nonlinear Optimization with Binary Constraints}

\author{Lianghai Xiao,~Yitian Qian and Shaohua Pan
\thanks{
The first two authors contributed equally to this work.

Lianghai Xiao is with College of Information Science and Technology, Jinan University, Guangzhou, China,
~Yitian Qian is with Department of Data Science and Artificial Intelligence, The Hong Kong Polytechnic University, Hong Kong, China,
and Shaohua Pan is with the School of Mathematics, South China University of Technology, Guangzhou, China.
(e-mail: xiaolh$@$jnu.edu.cn, yitian.qian@polyu.edu.hk, shhpan$@$scut.edu.cn. Corresponding author: Yitian Qian.) }
}
\markboth{Journal of \LaTeX\ Class Files,~Vol.~14, No.~8, August~2021}%
{Shell \MakeLowercase{\textit{Xiao et al.}}: A relaxation method for binary optimizations on constrained Stiefel manifold}


\maketitle

\begin{abstract}

We study binary optimization problems of the form 
\(
\min_{x\in\{-1,1\}^n} f(Ax-b)
\)
with possibly nonsmooth loss \(f\). Following the lifted rank-one semidefinite programming (SDP) approach\cite{qian2023matrix},  we develop a majorization-minimization algorithm by using the difference-of-convexity (DC) reformuation for the rank-one constraint and the Moreau envelop for the nonsmooth loss.
We provide global complexity guarantees for the proposed \textbf{D}ifference of \textbf{C}onvex \textbf{R}elaxation \textbf{A}lgorithm (DCRA)  and show that it produces an approximately feasible binary solution with an explicit bound on the optimality gap.
Numerical experiments on synthetic and real datasets confirm that our method achieves superior accuracy and scalability compared with existing approaches.

\end{abstract}


\section{Introduction}
\IEEEPARstart{C}{onsider} the following binary optimization problem:
\begin{equation}\label{eq:ori}
\underset{x\in \{-1, 1\}^n}{\min} f(Ax-b),
\end{equation}
where $f : \mathbb{R}^r \to \mathbb{R} $ is a  locally Lipschitz continuous (possibly nonsmooth) function, $A \in \mathbb{R}^{r \times n}$, and $b \in \mathbb{R}^r$. 
This model subsumes a wide range of important binary optimization problems across engineering, computer science, and applied mathematics. 
We highlight several representative instances.   

\textbf{Maximum 2-SAT.\cite{AspvallPlassTarjan1979}}  
Consider $m$ logical clauses, each involving two variables. The goal is to maximize the number of satisfied clauses, equivalently to minimize the number of violations:
\[
\underset{x \in \{-1,1\}^n}{\min}\; \sum_{j=1}^m \ell_j(x),
\]
where each $\ell_j(x)$ is a nonsmooth loss indicating whether clause $j$ is violated.  

\textbf{Binary compressed sensing and robust regression.} \cite{chevalier2024compressed, lotfi2020compressed}
In sparse signal recovery, one seeks binary signal that fits linear measurements:
\[
\underset{x \in \{0,1\}^n}{\min} \; \|Ax-b\|_1 + \rho \|x\|_1,
\]
where the $\ell_1$ norm promotes robustness and sparsity simultaneously, and a simple change of variables relates the $\{0,1\}^n$ and $\{-1,1\}^n$ encodings.  
A closely related special case is absolute value regression:
\[
\underset{x \in \{-1,1\}^n}{\min} \; \|Ax-b\|_1,
\]
which directly corresponds to \eqref{eq:ori} with $f(z) = \|z\|_1$.  

\textbf{MIMO detection in communications.}  \cite{nguyen2021linear, shao2024one}
For a received signal $y \in \mathbb{R}^r$ and channel matrix $H \in \mathbb{R}^{r\times n}$, symbol detection is modeled as
\[
\underset{x \in \{-1,1\}^n}{\min} \; \|Hx - y\|_2^2,
\]
which is a binary least-squares problem.  

\textbf{Hashing and binary embedding.} \cite{Xiong21, xiao2026relaxation, gui2016supervised}
In large-scale retrieval, supervised hashing learns binary codes by solving
\[
\underset{X \in \{-1,1\}^{n\times r}, \, W \in \mathbb{R}^{d\times r}}{\min} \;\; \|B - W X^\top\|_1 + \tfrac{\delta}{2}\|W\|_F^2,
\]
where $B$ is a data- or label-related matrix, $W$ is a linear projection, and $X$ is the binary code matrix. This fits the form~\eqref{eq:ori} after vectorization.

These examples illustrate that \eqref{eq:ori} captures both classical NP-hard combinatorial problems (QAP\cite{Burkard13}, 2-SAT\cite{AspvallPlassTarjan1979}) and modern applications in signal processing, statistics, communications, and machine learning. 
In particular, nonsmooth formulations such as $\ell_1$ regression and supervised hashing underscore the importance of developing scalable methods for \eqref{eq:ori} beyond the smooth or quadratic setting.


\subsection{Related Work}\label{1partB}
Despite its broad relevance, solving problems of the form~\eqref{eq:ori} remains a major challenge due to the interplay between combinatorial constraints and nonsmooth objectives. Classical approaches include exact and near-exact methods such as cutting-plane and branch-and-cut algorithms~\cite{Wu15,Glover10}. Cutting-plane methods iteratively strengthen relaxations by adding valid inequalities, while branch-and-cut combines branching with dynamically generated cuts to tighten linear or semidefinite relaxations. Although these schemes can provide strong bounds and, in principle, guarantee global optimality, they are often computationally prohibitive: in the worst case, they may require exploring up to $2^n$ subproblems, which limits their practicality for large-scale instances.

To improve scalability, a large body of work has focused on continuous relaxation methods, where the binary feasible set is replaced by a convex surrogate. Linear programming (LP) relaxation~\cite{komodakis2007approximate} substitutes $\{-1,1\}^n$ with the box constraints $-1 \le x \le 1$, followed by a rounding step to recover an integral solution. SDP relaxation~\cite{Goemans95} is more powerful: for the Max-Cut problem, an SDP relaxation combined with randomized rounding achieves an approximation ratio of at least $0.87856$. Subsequent works have strengthened SDP relaxations using cutting-plane inequalities~\cite{Helmberg98,Rendl10,Krislock14} and doubly nonnegative relaxations~\cite{Peng13,Fu18,qian2023matrix}. These approaches demonstrate the effectiveness of convex relaxations, although their practical performance can be sensitive to problem structure and to the quality of the rounding schemes.

Another important line of research uses penalty-based methods to enforce binary constraints. Early approaches studied smooth penalty formulations~\cite{Anjos02,Burer01,Chardaire94}. Later, strategies based on mathematical programs with equilibrium constraints (MPECs)~\cite{bi2014exact,Yuan17} were proposed and applied to tasks such as binary hashing and subgraph discovery. More recently, augmented Lagrangian and penalty decomposition frameworks have been developed for general binary optimization~\cite{Xiong21}.

Despite these advances, a general-purpose framework for nonsmooth binary optimization is still relatively limited. Many existing methods are tailored to specific problem classes, especially those with smooth or quadratic objectives. Penalty methods~\cite{Xiong21,Yuan17} often assume differentiability to establish convergence, while convex relaxations such as SDP and doubly nonnegative formulations are most effective for quadratic objectives where the structure can be fully exploited. When the objective $f$ is nonsmooth and highly sensitive to perturbations in $Ax$, relaxation-and-rounding strategies may lead to a significant deterioration in solution quality. This motivates the development of new approaches that combine the strengths of SDP relaxations with exact penalty techniques, enabling efficient handling of nonsmooth objectives while maintaining scalability and theoretical guarantees.
\subsection{Contributions}   
 
The contributions of this paper are threefold. 
First, we extend the rank-one DC penalty framework of~\cite{qian2023matrix} to problems of the form~\eqref{eq:ori} with general (possibly nonsmooth) objectives, thereby covering a broader range of applications such as robust regression, binary compressed sensing, and supervised hashing. 
We observe that $x \in \{-1,1\}^n$ if and only if the matrix
$$
X = \begin{bmatrix}1 & x^\top \\ x & xx^\top \end{bmatrix}
$$
is positive semidefinite, rank-one, and satisfies $\text{diag}(X) = e$. 
This observation allows us to reformulate~\eqref{eq:ori} as
\begin{equation}\label{eq:eqv1}
	\underset{X \in \mathbb{S}^{p}_+}{\min} \left\{ f\lf({\cal A}(X) - b\rg) \; \text{s.t.} \; \text{rank}(X) \leq 1, \; \text{diag}(X) = e \right\},
\end{equation}
where $p := n+1$ and ${\cal A} : \mathbb{S}^p \to \mathbb{R}^r$ is the linear operator defined by
\(
  {\cal A}(Z) = \tfrac{1}{2}A(z_2 + z_3)
\)
for
\(
  Z = \begin{bmatrix}z_1 & z_2^\top \\ z_3 & z_4 \end{bmatrix},
\)
with \( z_1 \in \mathbb{R} \), \( z_2, z_3 \in \mathbb{R}^n \), and \( z_4 \in \mathbb{R}^{n \times n} \).
We further encode the rank-one constraint using the DC representation $\|X\|_* - \|X\| = 0$, leading to a DC-constrained SDP reformulation 
\begin{equation}\label{eq:eqv2}
	\underset{X\in\mathbb{S}^{p}_+}{\min}\lf\{f\lf({\cal A}(X)-b\rg),\; {\rm s.t.} \;\|X\|_* - \|X\| = 0,\,{\rm diag}(X) = e\rg\}.
\end{equation}
 
Second, we develop a majorization–minimization scheme (DCRA) that leverages the Moreau envelope of $f$,  which admits closed-form updates in subproblems and avoids repeated eigen-decompositions. 

Third, we show that the inner iteration scheme has $O(1/\varepsilon^{2})$ worst-case complexity in terms of a stationarity residual, and the overall penalty algorithm terminates in finitely many steps. 
At any normal termination point, a simple rank-one projection of the lifted iterate yields an approximately feasible binary solution together with an explicit, computable bound on its optimality gap.

\subsection{Organization}

The remainder of this paper is structured as follows. 
Section \ref{sec:pre} sets notation and preliminaries. 
Section~\ref{sec:equiv} develops the DC-penalized model and its Burer–Monteiro factorization, and clarifies its equivalence to the lifted rank-one SDP.
Section~\ref{sec:algo} presents the smoothed MM algorithm on the unit-column sphere, including the closed-form inner update. 
Section~\ref{sec:conv} provides the convergence and complexity analysis of the proposed framework, as well as a rank-one projection guarantee with an explicit lower-bound gap.
Section~\ref{sec:exp} reports numerical results, including parameter sensitivity studies and experiments on synthetic problems, binary compressed sensing, and supervised hashing.
Finally, Section~\ref{sec:conc} concludes the paper.

\section{Preliminaries}\label{sec:pre}
\IEEEPARstart{T}{hroughout} this paper, we use the following notation and definitions. Let $\mathbb{R}^{n_1 \times n_2}$ denote the space of all real $n_1 \times n_2$ matrices, equipped with the trace inner product $\langle X, Y \rangle := \operatorname{Tr}(X^\top Y)$ and its associated Frobenius norm $\|X\|_F := \sqrt{\langle X, X \rangle}$. Let $\mathbb{O}^{p}$ denote the set of all $p \times p$ orthonormal matrices. For a positive integer $k$, we write $[k] := \{1, 2, \ldots, k\}$.
For a matrix $X \in \mathbb{R}^{n_1 \times n_2}$, we denote its spectral norm and nuclear norm by $\|X\|$ and $\|X\|_*$, respectively. For an index set $J \subseteq [n_2]$, $X_J$ denotes the submatrix of $X$ consisting of columns $X_j$ with $j \in J$. The closed ball centered at $X$ with radius $\varepsilon > 0$ under the Frobenius norm is denoted $\mathbb{B}(X, \varepsilon)$.

Denoted by $\mathbb{S}^p$ the set of symmetric matrix. 
For a symmetric matrix $X \in \mathbb{S}^p$, we use $\operatorname{diag}(X) \in \mathbb{R}^p$ to denote the vector of its diagonal entries, for any vector $z \in \mathbb{R}^p$, $\operatorname{Diag}(z) \in \mathbb{S}^p$ denotes the diagonal matrix with $z$ on its diagonal. We denote the eigenvalues of $X$ by $\lambda(X) := (\lambda_1(X), \ldots, \lambda_p(X))$, arranged in nonincreasing order, and , and $\Lambda_m(X ) := \operatorname{Diag}(\lambda_1(X ), \ldots, \lambda_m(X ))$ the diagonal matrix constructed by the first $m$ eigenvalues of $X$. Define the set of orthogonal diagonalizers by
$
\mathbb{O}(X) := \left\{ P \in \mathbb{O}^p \mid X = P \operatorname{Diag}(\lambda(X)) P^\top \right\}. 
$ 
For $P \in \mathbb{O}(X)$, we use $P_I \in \mathbb{R}^{p \times m}$ to denote the submatrix of $P$ containing its first $m$ columns. 
The identity matrix, the all-ones vector, and the all-ones matrix are denoted by $I$, $e$, and $E$, respectively, with their dimensions clear from context.

Given a closed set $\Delta \subseteq \mathbb{R}^{n_1 \times n_2}$, the indicator function $\delta_\Delta$ is defined by
$$
\delta_\Delta(z) := \begin{cases}
0, & \text{if } z \in \Delta, \\
+\infty, & \text{otherwise}.
\end{cases}
$$ 
The distance from a point $X \in \mathbb{R}^{n_1 \times n_2}$ to $\Delta$ is denoted $\operatorname{dist}(X, \Delta)$, and is measured in Frobenius norm.
For notational convenience, we define:
{\small
\begin{equation}\label{psi-fun}
 \widetilde{\psi}(V) := \psi(V^\top V),\;\text{for } V \in \mathbb{R}^{m \times p}, \text{ with } \psi(X) := -\|X\|.
\end{equation}}
We now recall the definitions of the regular and (limiting) subdifferentials from \cite{RW98}, which will be used to analyze nonsmooth functions. 

 \begin{definition}\label{Gsubdiff-def}
Let $h: \mathbb{R}^n \to (-\infty, +\infty]$ and $x \in \mathbb{R}^n$ be such that $h(x) < +\infty$. The regular subdifferential of $h$ at $x$, denoted $\widehat{\partial} h(x)$, is defined as
{\footnotesize
$$
\widehat{\partial} h(x) := \left\{ v \in \mathbb{R}^n \ \middle| \
\liminf_{x' \to x,\, x' \neq x} \frac{h(x') - h(x) - \langle v, x' - x \rangle}{\|x' - x\|} \ge 0 \right\}.
$$
 }
The (limiting) subdifferential of $h$ at $x$ is given by
\[
\partial h(x) := \Bigl\{ v \in \mathbb{R}^n \ \Bigm|\ 
\begin{aligned}
& \exists\, x^k \to x \text{ with } h(x^k) \to h(x), \\
& v^k \in \widehat{\partial} h(x^k),\; v^k \to v \ {\rm as}\ k\to\infty
\end{aligned}
\Bigr\}.
\]
\end{definition}
It follows from the definition that $\widehat{\partial} h(x) \subseteq \partial h(x)$, and both sets are closed; $\widehat{\partial} h(x)$ is also convex. When $h$ is convex, both subdifferentials reduce to the usual convex subdifferential. A point $\bar{x}$ is called a critical point of $h$ if $0 \in \partial h(\bar{x})$; we denote the set of all such points by $\operatorname{crit}\, h$.
For an indicator function $h = \delta_\Delta$, the subdifferentials reduce to the regular normal cone $\widehat{\mathcal{N}}_\Delta(x)$ and the normal cone $\mathcal{N}_\Delta(x)$, respectively.

\section{Reformulation}\label{sec:equiv}
\IEEEPARstart{T}{he} feasible set of problem \eqref{eq:eqv2} remains combinatorial, as it consists rank-one positive semidefinite (PSD) matrices with binary structure. Notably, the rank constraint  can be written as a DC constraint $\|X\|_* - \|X\| = 0$, but handling a DC constraint is generally more difficult than handling a DC objective. To address this difficulty, we consider a penalty reformulation \eqref{eq:eqv2} that incorporates the DC constraint into the objective function 
\begin{equation}\label{eq:pen}
\underset{X \in \Omega}{\min}  f\big({\cal A}(X) - b\big) + \rho \big( \|X\|_* - \|X\| \big),
\end{equation}
where $\rho > 0$ is a penalty parameter, and
\(\Omega := \{ X \in \mathbb{S}_+^p \mid \operatorname{diag}(X) = e \}.\)

A key fact is that the penalized model \eqref{eq:pen} is an exact-penalty reformulation of \eqref{eq:eqv2}. That is, there exists  \(\rho^* = (1+2p)L_f\)  such that, for all \(\rho\ge\rho^*\), the global minimizers of \eqref{eq:pen} coincide with those of \eqref{eq:eqv2} by\cite[Proposition 2.3 and Theorem 3.1]{BiPan16}, and hence with the solutions of the rank-one SDP \eqref{eq:eqv1}. Moreover, for any $\rho>\rho^*$, every feasible point  $X\in\mathcal F$ of   \eqref{eq:eqv2} is not only locally optimal for \eqref{eq:eqv2} but also  a strict local minimizer of  \eqref{eq:pen} (see  \cite[Prop.~1(ii)]{qian2023matrix}). 

In practice, however, solving problem \eqref{eq:pen} is nontrivial. Even if the threshold $\rho^*$ is known, one typically cannot expect to obtain global optima for a fixed $\rho \ge \rho^*$. Moreover, existing convex relaxation methods require full eigenvalue decompositions at each iteration. This causes a major computational bottleneck that impedes scalability to large-scale settings.
To overcome this, we adopt a matrix factorization paradigm \cite{Burer01, Burer03}, which has gained renewed interest in low-rank optimization and semidefinite relaxation (e.g., \cite{SunLuo16, Li18}). Specifically, we introduce the following penalized factorized formulation
{\small
\begin{equation}\label{eq:eqv3}
\underset{V \in \mathcal{S}}{\min} \;\Phi(V):= f\big({\cal A}(V^\top V) - b\big) 
+ \rho \left( \|V\|_F^2 - \|V\tp V\| \right),
\end{equation}}
where $\mathcal{S} := \left\{ V \in \mathbb{R}^{m \times p} \mid \|V_j\| = 1 \text{ for } j = 1, \ldots, p \right\}$, $2 \le m < p$, and $V_j$ denotes the $j$th column of $V$. The penalty term $\|V\|_F^2 - \|V\|^2$ corresponds to the surrogate of the original rank-one regularizer, and the column-wise normalization ensures  $\operatorname{diag}(X) = e$. We require $m \ge 2$ since the penalty term vanishes when $m = 1$.
Although $\|V\|_F^2 = p$ for all $V \in \mathcal{S}$, we retain this term to explicitly encode the penalization for the rank surrogate $(\|X\|_* - \|X\|)$. 

By  \cite[Sec.~3]{qian2023matrix}, if $V^* \in \mathcal{S}$ is globally optimal for \eqref{eq:eqv3}, then $X^* := V^{*\top} V^*$ is globally optimal for \eqref{eq:pen} with $\operatorname{rank}(X^*) \le m$. 
Conversely, given any global minimizer $X^*$ of \eqref{eq:pen} with rank $r^*$, one can obtain a global optima  to \eqref{eq:eqv3} by letting $V^* := \sqrt{\Lambda_m(X^*)} P_I^{* \top}$. 
For convenience, we refer to problem \eqref{eq:eqv3} as a \textit{DC penalized matrix program}, even though it is not a true DC program due to the nonconvexity of the feasible set $\mathcal{S}$.  
 
\section{Algorithmic Framework}\label{sec:algo}

\IEEEPARstart{T}{o}  address the nonsmoothness of the objective function in \eqref{eq:eqv2}, we consider its Moreau envelope smoothing. For a  proper, lower semi-continuous function $f:\mathbb{R}^q \to (-\infty,+\infty]$, its Moreau envelope of parameter $\delta > 0$ and the corresponding proximal mapping are defined as
\begin{equation}\label{eq:moreau-env}
{\rm env}_{\delta f}(x) := \underset{y\in\mathbb{R}^q}{\min} \left\{ f(y) + \frac{1}{2\delta} \|x - y\|^2 \right\}, \quad {\rm and},
\end{equation} 
\begin{equation}\label{eq:prox-general}
{\rm prox}_{\delta f}(x) := \underset{y\in\mathbb{R}^q}{\arg\min} \left\{ f(y) + \frac{1}{2\delta} \|x - y\|^2 \right\}, 
\end{equation}
respectively. If $f$ is convex, then ${\rm env}_{\delta f}$ is continuously differentiable, and its gradient is given by
\(
\nabla\, {\rm env}_{\delta f}(x) = \frac{1}{\delta} \lf(x - {\rm prox}_{\delta f}(x)\rg).
\)
Using this smoothing technique, we consider the following smoothed problem:
{\footnotesize
\begin{equation}\label{eq:smooth-fact}
 \underset{V \in \mathcal{S}}{\min}  \Phi_\delta(V):=\widetilde{f}(V) + \rho \left( \|V\|_F^2 + \widetilde{\psi}(V)  \right) ,
\end{equation}}
where   $\widetilde{f}(V) := {\rm env}_{\delta f}(\mathcal{A}(V^\top V) - b)$, and $\widetilde{\psi}(V)$ is defined in \eqref{psi-fun}.
 We propose the following smoothed relaxation approach by seeking a finite number of stationary points \eqref{eq:smooth-fact} associated to an increasing $\rho$.

\begin{algorithm}[H]
\caption{\textbf{D}-\textbf{C} \textbf{R}elaxation \textbf{A}lgorithm (DCRA) based on model \eqref{eq:smooth-fact}}
\label{alg:dc_pen}
\begin{algorithmic}[1] 
\State Initialization: select an appropriate integer $m\in [2, p]$, a small $\epsilon\in(0, 1)$, and an appropriately large $k_{\max}\in {\mathbb N}$. Choose $\rho_{\max}>0, \rho_0>0, \sigma>1$,  \( \delta > 0 \), and a starting point $V^0\in {\cal S}$ with full rank
. 
        
\For{$k=0,1,2,\ldots, k_{\max}$}
\State\label{Astep3} From $V^k\in{\cal S}$, seek an approximate stationary point $V^{k+1}$ of the following problem:
 
\State {\footnotesize \begin{equation}\label{eq:alg1}
\underset{V \in \mathcal{S}}{\min} \Phi_k (V) :=\widetilde{f}(V) + \rho_k \left( \|V\|_F^2 + \widetilde{\psi}(V)  \right) ,
\end{equation}}
\State If $\|V^{k+1}\|_F^2 - \|V^{k+1}\|^2 \le \epsilon$, then stop. 
\State Otherwise, let $\rho_{k+1}\leftarrow\min\{\sigma\rho_k, \rho_{\max}\}$.\label{step:rho}
\EndFor  
\end{algorithmic}
\end{algorithm}
We next consider the subproblem \eqref{eq:alg1}.
\begin{definition}\label{def:DC-stationary}
Given $\epsilon > 0$, a point $V^\star\in\mathcal S$ is an $\epsilon$-stationary point
of \eqref{eq:alg1} if there exists $\Gamma^\star\in\partial\widetilde\psi(U)$ for some $U\in {\cal S}\cap\mathbb{B}(V, \epsilon)$ such that
{\small
\[
{\cal G}_k(V^\star):=\operatorname{dist}\!\Big(0,\ \nabla\widetilde f(V^\star) + 2\rho_k V^\star + \rho_k\Gamma^\star+ N_{\mathcal S}(V^\star)\Big) \le \epsilon.
\]}
\end{definition}
To compute a stationary point of \eqref{eq:alg1}, we linearize $\widetilde{f}$ and $\widetilde{\psi}$ at the current iterate $V^{k,l}$, and obtain the subproblem:
{\footnotesize
\begin{align}\label{eq:lin-sub}
\underset{V \in \mathcal{S}}{\min} \;
& 
\widetilde\Phi_k (V;V^{k,l}):=\widetilde{f}(V^{k,l})
  + \langle \nabla \widetilde{f}(V^{k,l}), V - V^{k,l} \rangle \notag\\
& 
+ \rho_k \left[ \|V\|_F^2
  + \widetilde{\psi}(V^{k,l})
  + \langle \Gamma^{k,l}, V - V^{k,l} \rangle \right] \notag\\
& + \frac{L}{2} \|V - V^{k,l}\|_F^2.
\end{align}
}
where $\Gamma^{k,l} \in \partial \widetilde{\psi}(V^{k,l})$ and $L \ge L_{\widetilde{f}}$. 
 This subproblem admits a closed-form solution:
{\footnotesize
\begin{equation}\label{eq:lin-sol}
V^{k,l+1} := \operatorname{proj}_{\mathcal{S}} \left(
\frac{1}{2\rho_k + L} \left( L V^{k,l} - \nabla \widetilde{f}(V^{k,l}) - \rho_k \Gamma^{k,l} \right)
\right),
\end{equation}}
where $\operatorname{proj}_{\mathcal{S}}(U):= {\rm normalize}(U)$ projects column-wise onto the unit sphere.

\begin{algorithm}[H]
\caption{Solver for the Subproblem \eqref{eq:alg1}}\label{AlgA}
\begin{algorithmic}[1] 
\State \textbf{Input:} current outer iterate \( V^{k,0} := V^{k} \in \mathcal{S} \), \( \rho_k > 0 \), tolerance \( \epsilon_v > 0 \), and  \( l_{\max} \)\label{alg:step1}
\For{$l = 0, 1, 2, \ldots, l_{\max}$}
    \State Compute \( \Gamma^{k,l} \in \partial \widetilde{\psi}(V^{k,l}) \) .
    \State Compute:
    \State  $ \nabla \widetilde{f}(V^{k,l}) = 2 V^{k,l} \mathcal{A}^* \nabla {\rm env}_{\delta f}(\mathcal{A}((V^{k,l})^\top V^{k,l}) - b) $.
    \State Update $V^{k,l+1}$ via \eqref{eq:lin-sol}. 
    \If{ \( \|V^{k,l+1} - V^{k,l}\|_F \le \epsilon_v \), } \textbf{break}
    \EndIf
\EndFor
\State \textbf{Output:} Approximate solution \( V^{k+1} := V^{k,l+1} \).
\end{algorithmic}
\end{algorithm}
 
By \cite[Lemma 1, Sec. 4]{qian2023matrix}, we may take
$$
\Gamma^{k,l} = \Gamma(V^{k,l}):=-2 V^{k,l} P_1(V^{k,l}) P_1(V^{k,l})^\top,
$$
where $P_1(V)$ is the leading eigenvector of $V^\top V$, obtainable from a partial SVD, whose computation cost is cheap since $V^{k,l}$ has less rows. 
\vspace{0.5cm}

\section{Convergence Analysis}\label{sec:conv}
\IEEEPARstart{I}{n} this section, we will provide the convergence result for our method. Define the residual
\[
  \mathcal G(V, \Gamma)
  :=\operatorname{dist}\!\big(0,\ \nabla\tilde f(V)+2\rho_k V+\rho_k\Gamma+N_{\mathcal S}(V)\big).
\]  
In our algorithmic context we will use $\Gamma=\Gamma^{k,l}\in\partial\tilde\psi(V^{k,l})$ and evaluate the residual at $V^{k,l+1}$.

\begin{proposition}\label{prop:inner-descent}
Let $\{V^{k,l}\}_{l\in\mathbb{N}}\subset\mathcal{S}$ be the sequence generated by Algorithm~\ref{AlgA}
for a fixed $k\in \mathbb{N}$. Then the following statements hold:
\begin{enumerate} 
\item[(i)]
For all $l\ge0$,
{\small
\begin{equation}\label{eq:nonincr}
  \Phi_k(V^{k,l})-\Phi_k(V^{k,l+1})
  \ \ge\
  \tfrac{L-L_{\tilde f}}{2}\,\|V^{k,l+1}-V^{k,l}\|_F^2.
\end{equation} }

\item[(ii)]
The sequence $\{V^{k,l}\}\subset\mathcal S$ is bounded, and $\{\Phi_k(V^{k,l})\}$ is nonincreasing and convergent.

\item[(iii)]  For all $l\in \mathbb{N}$, there exists $C>0$ such that
{\small
\begin{equation}\label{eq:prop4}
\mathcal G\big(V^{k,l+1}, \Gamma^{k,l}\big)
\le C \|V^{k,l+1}-V^{k,l}\|_F.
\end{equation}}
\end{enumerate}
\end{proposition}

\begin{theorem}\label{thm:inner-complexity}
Let  
\(
  \Phi_{\delta}^\star := \min_{V\in\mathcal S}\,  \Phi_{\delta}(V)
\)
be the optimal value of  the smoothed problem \eqref{eq:smooth-fact}.
Fix outer iteration $k\in\mathbb{N}$. For any $l\in [T]$ with $T\in\mathbb N$, the inner iterates generated by Algorithm~\ref{AlgA} satisfy
{\footnotesize
\begin{equation}\label{eq:itr_cmplx}
  \underset{0\le l<T}{\min}\; \mathcal G\big(V^{k,l+1}, \Gamma^{k,l}\big)
  \ \le\
  \frac{C}{\sqrt{T}}\,
  \sqrt{\frac{2\big(\Phi_k(V^{k})-\Phi_{\delta}^\star \big)}{L-L_{\tilde f}}}.
\end{equation}
}
In particular, if
\(
T  \ge \left\lceil \frac{2C^2}{L-L_{\tilde f}} 
\frac{\Phi_k(V^{k})-\Phi_{\delta}^\star}{\epsilon^2}\right\rceil,
\)
then the inner loop attains an $\epsilon$-stationary point within $O(\epsilon^{-2})$ iterations, and hence Algorithm \ref{AlgA} terminates after finitely many inner iterations.
\end{theorem}

Next, we show that Algorithm~\ref{alg:dc_pen} terminates normally. Lemma~\ref{lem:complexity} ensures that the gap $\|V^{k,1}\|_F^2 - \|V^{k,1}\|$ is sufficiently small under suitable conditions, and Proposition~\ref{prop:inner-progress} then guarantees that Algorithm~\ref{alg:dc_pen} exits properly. Although the analysis in this part follows the techniques of~\cite{qian2023matrix}, the results we obtained are different.

\begin{lemma}\label{lem:complexity}
Fix $k\in\mathbb{N}$. Suppose that there exists 
$P\in\mathbb{O}\!\big((V^{k})^{\top}V^{k}\big)$ such that 
its leading eigenvector $P_1$ has no zero entries. 
Define 
\(
\mu\ :=\ \min_{1\le j\le p}\,|P_{j1}|\ >0.
\)
Set
\(
\varpi\ :=\ (L+2L_{\widetilde f})\sqrt{p}
\ +\ \big\|\nabla\widetilde f\!\big(\tfrac{1}{\sqrt m}E\big)\big\|_F.
\)
Take
\(
\rho_k\ \ge\ \overline\rho\ :=\ \frac{Lp+\varpi}{\mu\,c_0}, 
\) 
with $c_0\in(0,1)$, 
then the first inner iterate $V^{k,1}$ generated by Algorithm~\ref{AlgA} satisfies
\[
\|V^{k,1}\|_F^2-\|V^{k,1}\|^2
\ \le\ c_0^2.
\]

\end{lemma}

\begin{proposition}\label{prop:inner-progress}
Fix $k\in\mathbb{N}$ and let $\{V^{k,l}\}_{l\ge0}\subset\mathcal S$ be the inner sequence produced by Algorithm~\ref{AlgA} with curvature $L\ge L_{\widetilde f}$ and $V^{k,0}=V^{k}$.
Let $\varepsilon\in(0,c_0^2]$ be a tolerance, where $c_0\in(0,1)$ is the same constant as in Lemma~\ref{lem:complexity}.
Define
\(
\varpi_{\mathrm{lin}}
:= 4L_{\widetilde f}p  +  2\sqrt p\,\Big\|\nabla\widetilde f\!\Big(\tfrac{1}{\sqrt m}E\Big)\Big\|_F
 + \tfrac{5}{2} L p
\), 
\(
\eta(c_0,p):=2\Big(1-\tfrac{1}{p}\Big)\sqrt{p-c_0^2} + 2\sqrt{1-c_0^2},
\)
and set
\(
\widehat{\rho}\ :=\ \max\lf\{\frac{2\,\varpi_{\mathrm{lin}}}{\eta(c_0,p)}, \frac{2\varpi_{\mathrm{lin}}}{\varepsilon}\rg\}.
\)
Then the following hold:
\begin{itemize}
\item[(i)] If $\rho_k\ge \widehat{\rho}$, then for every $l\in\mathbb N$ such that
$\frac{\varepsilon}{2}\le p-\|V^{k,l}\|^2\le c_0^2$, one has
\begin{equation}\label{eq:gap_v}
\|V^{k,l+1}\|^2
\ \ge\
\|V^{k,l}\|^2\ +\ \tfrac{1}{2}\,\eta(c_0,p).
\end{equation}

\item[(ii)] Suppose there exists $P\in\mathbb O\big((V^{k,0})^\top V^{k,0}\big)$ whose leading eigenvector $P_1$ has no zero entries, and if
$\rho_k>\max\{\widehat{\rho},\overline\rho\}$ with $\overline\rho :=  \frac{Lp+\varpi}{\mu\,c_0}$ from Lemma \eqref{lem:complexity}, then there exists
\(
1\ \le\ \overline l\ \le\ \Big\lceil\frac{2c_0^2-\varepsilon}{\eta(c_0,p)}\Big\rceil+1
\) 
such that for all $l\ge \overline l$,
\[
\|V^{k,l}\|_F^2-\|V^{k,l}\|^2\ =\ p-\|V^{k,l}\|^2\ \le\ \varepsilon.
\]

\end{itemize}
\end{proposition}

The following corollary is an easy conclusion from Theorem \ref{thm:inner-complexity} and Proposition \ref{prop:inner-progress}. 

\begin{corollary}
\label{thm:overall-complexity}
Suppose there exists $P\in\mathbb O\big((V^{k,0})^\top V^{k,0}\big)$ whose leading eigenvector $P_1$ has no zero entries, and $\varepsilon_{v}$ is given in Algorithm \ref{AlgA}. The total number of iterations performed by Algorithm~\ref{alg:dc_pen} is bounded by
\(
\mathcal O\!\Big(\tfrac{\log(\rho_\star/\rho_0)}{\varepsilon_{v}^{2}}\Big) 
\)
with $ \rho_\star := \max\{\widehat\rho,\overline\rho\}$, where $\widehat\rho$ and $\overline\rho$ are given in Proposition~\ref{prop:inner-progress}
and Lemma~\ref{lem:complexity}, respectively.
\end{corollary}

Theorem \ref{thm:feas-obj} states that the rank-one projection of the normal output of Algorithm \ref{alg:dc_pen} is approximately feasible and near-optimal. 

\begin{theorem}\label{thm:feas-obj}
Let $\overline{k}$ be an outer index at which Algorithm \ref{alg:dc_pen} terminates normally, in which the subproblem is solved by Algorithm \ref{AlgA}. 
Let $X^{\overline{k}}=(V^{\overline{k}})^\top V^{\overline{k}}$ and let $P\in\mathbb O(X^{\overline{k}})$ diagonalize $X^{\overline{k}}$.
Define the rank-one projection $x^{\overline{k}}:=\|V^{\overline{k}}\|\,P_1$, where $P_1$ is the leading eigenvector of $X^{\overline{k}}$. Write 
\(
\hat{f}(X):= \mathrm{env}_{\delta f}(\mathcal A(X)-b). 
\) 
Assume $\hat{f}$ is Lipschitz on $\Omega$ with modulus $L_{\hat f}>0$.
Then, for any $k^*\in\{0,1,\ldots,\overline{k}-1\}$, the following hold:
\begin{align}
\|x^{\overline{k}} \circ x^{\overline{k}}-e\| & \le\ \varepsilon, \label{eq:feas-dcra}\\[-1mm]
\hat{f}\!\big(x^{\overline{k}}(x^{\overline{k}})^\top\big)\ -& \widetilde f(V^{k^*})
\ \le\ \rho_{\overline{k}}\|V^{\overline{k}}\|^2\ -\ \rho_{k^*}\,\frac{p}{r^*}
\notag\\[-2mm]
+& \sum_{j=k^*}^{\overline{k}-1}\!\big(\rho_j-\rho_{j+1}\big)\,\|V^{j+1}\|^2
\ +\ L_{\hat f}\,\varepsilon, \label{eq:obj-gap-dcra}
\end{align}
where $r^*:=\mathrm{rank}(V^{k^*})$. 
In particular, if $\widetilde f(V^{k^*})=\hat{f}\!\big((V^{k^*})^\top V^{k^*}\big)\le \upsilon^*$ 
for the optimal value $\upsilon^*$ of the original problem, then
\[
\hat{f} \big(x^{\overline{k}}(x^{\overline{k}})^\top\big) - \upsilon^*
\le \rho_{\overline k}(p-1) 
+ \rho_{k^*}\Big(1-\frac{p}{r^*}\Big) 
+ (L_{\hat f}-\rho_{\overline k})\epsilon.
\]
\end{theorem} 

\section{Numerical Experiment}\label{sec:exp}
\IEEEPARstart{W}{e} assess the performance of the proposed algorithm on both synthetic and real-world datasets, and benchmark it against several baseline methods. All experiments are implemented in MATLAB and executed on a workstation running 64-bit Windows, equipped with an Intel(R) Core(TM) i9-12900 CPU @ 2.40 GHz and 64 GB RAM.
The stopping criterion for all algorithms is either (i) reaching a maximum number of iterations or (ii) meeting a predefined tolerance for the change in objective value, whichever occurs first. Unless otherwise stated, algorithm parameters are set to default values.


\subsection{Parameter sensitivity analysis}

As the performance of Algorithm~\ref{alg:dc_pen} may be influenced by
(i) the number of columns $m$ in the initialization matrix $V_0$,
(ii) the initial penalty parameter $\rho$, and
(iii) the smoothing parameter $\delta$ in the subproblem formulation,
we conduct experiments to examine the sensitivity of its performance to these key parameters.
For each parameter, we consider two problem sizes: $(n, r) = (500, 300)$ and $(n, r) = (1000, 600)$. Results are averaged over 50 independent random trials.

For $m$, we vary its value over a logarithmic scale from $10^0$ to $10^3$.
For both problem sizes, the objective value improves rapidly as $m$ increases from small values, then stabilizes once $m$ exceeds roughly $1$–$10$.
When $m > 10$, the objective value becomes less stable for the larger problem $(n, r) = (1000, 600)$.
The runtime grows slowly for small $m$ but rises sharply for large $m$.

For $\rho$, moderate values yield stable, high-quality solutions with relatively low computation time.
Very small $\rho$ leads to unstable objective values and longer runtimes, while excessively large $\rho$ can also degrade performance.

For $\delta$, performance remains stable over a broad range of values.
Very large $\delta$ can reduce computation time, but often at the expense of solution quality.

\begin{figure}[h]
\centering
\subfloat[Effect of $m$]{
\includegraphics[width=\linewidth]{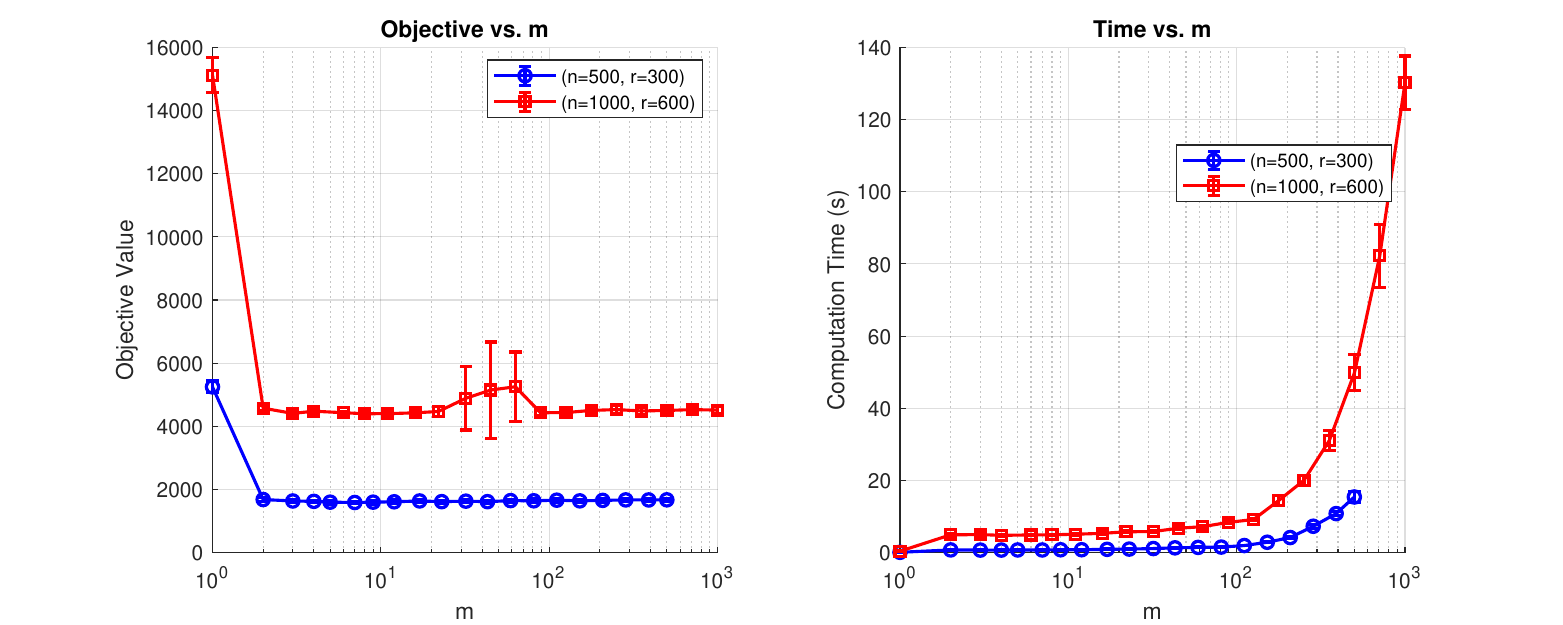}
\label{fig:sen_m}
} 

\vspace{-2mm}
\subfloat[Effect of the initial penalty parameter $\rho$]{
\includegraphics[width=\linewidth]{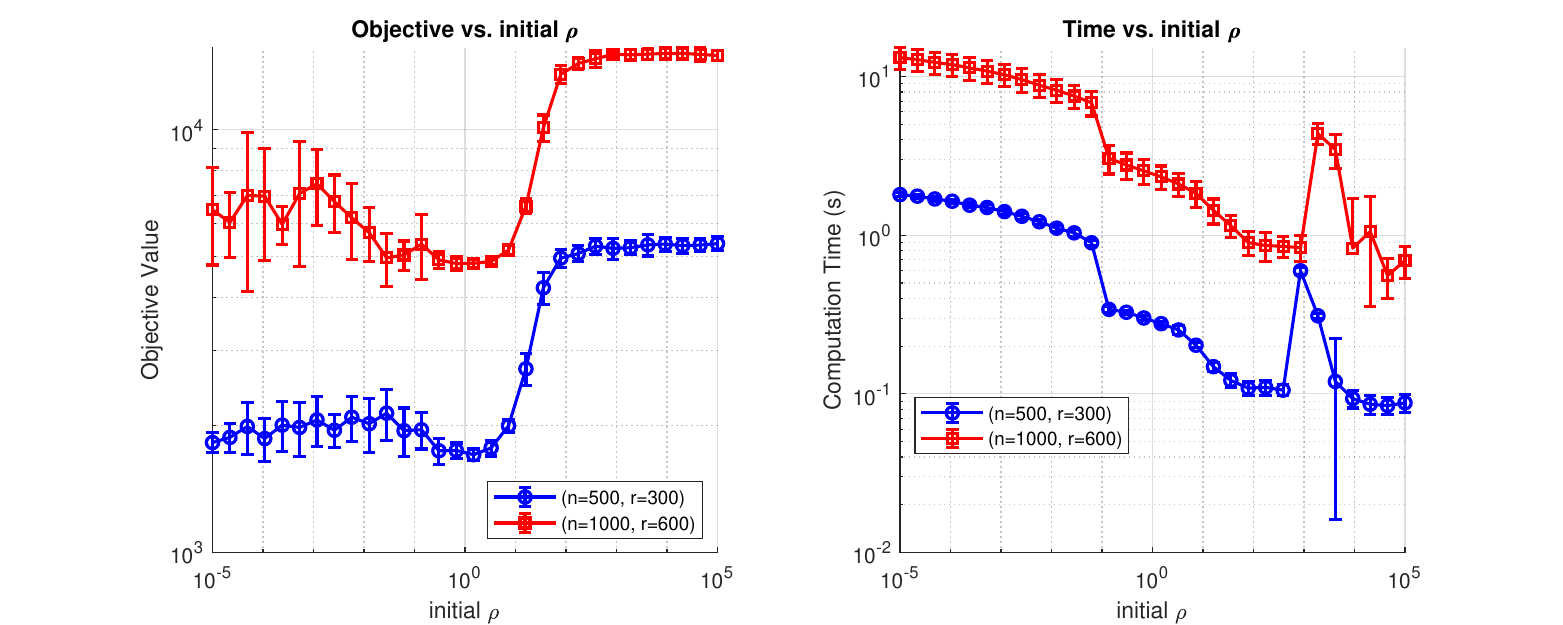}
\label{fig:sen_rho}
} 
 
\vspace{-2mm}
\subfloat[Effect of the smoothing parameter $\delta$]{
\includegraphics[width=\linewidth]{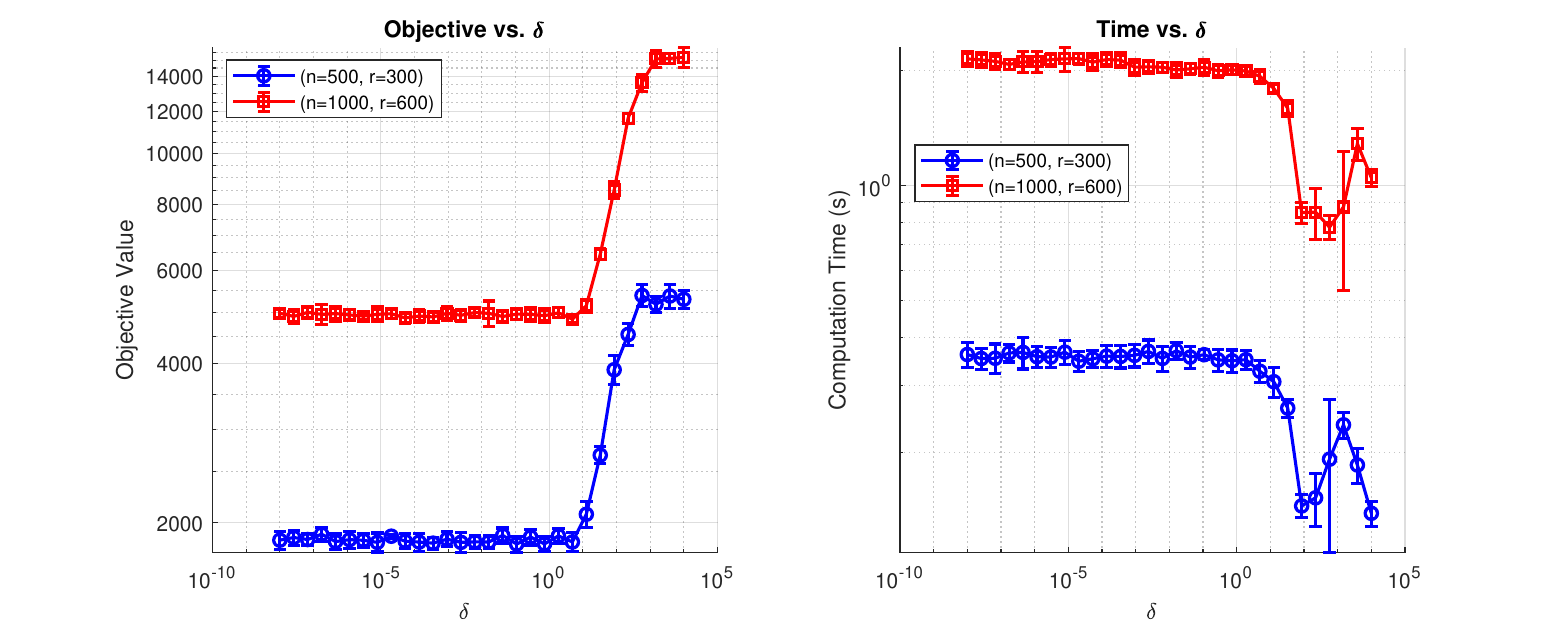}
\label{fig:sen_delta}
}
\caption{Parameter sensitivity analysis of the proposed algorithm for two problem sizes $(n,r)=(500,300)$ and $(n,r)=(1000,600)$. For each parameter setting, results are averaged over 10 trials, with error bars indicating one standard deviation. Left panels: objective value; Right panels: runtime. Both axes use logarithmic scales where appropriate.}
\label{fig:sen_all}
\end{figure}

\subsection{Performance on randomly generated problems}

We first evaluate the proposed algorithm on randomly generated instances to assess its performance and robustness. The objective function is
\begin{equation}\label{eq:oril1}
\min_{x} \|A x - b\|_1,
\end{equation}
where $A \in \mathbb{R}^{n \times r}$ and $b \in \mathbb{R}^n$.
For each problem size, we randomly generate $100$ independent instances, with each entry of $A$ and $b$ drawn from the standard normal distribution.
We test problem sizes $n\in [100, 3000]$ and $d\in[50, 2000]$. 
The parameters for DCRA are fixed as $\sigma = 1.2$, $m = 5$, and $\rho = 1$.

We compare DCRA with three baselines: the Gurobi solver, the MPEC-EPM method\cite{Yuan17}, and the L2-box ADMM method\cite{wu2018ell}. 
Performance is evaluated using the following four indicators:
\begin{itemize}
\item Win rate of DCRA: the proportion of instances in which the given method achieves a better objective value than the others;

\item Mean relative difference: the average value of $(\text{Obj}_{\text{DCRA}} - \text{Obj}_{\cdot})/\text{Obj}_{\cdot}$, where negative values indicate that DCRA achieves a smaller objective value;

\item Average objective value: the mean of $\|A x - b\|_1$ ; and

\item Average CPU time: the mean computation time over $100$ runs.

\end{itemize}

To apply Gurobi to our binary optimization model \eqref{eq:oril1}, we reformulate it as an equivalent mixed-integer linear programming problem(MILP):
{\small
$$
\begin{array}{rc}
\min_{x, z} & e\tp z \\
{\rm s.t.} & z \ge -(A(2x-1) - b), \\
           & z \ge A(2x-1) - b, \\
           & x \in \{0, 1\}^n, z\in \R^r.
\end{array}
$$}
This formulation allows us to exploit Gurobi’s branch-and-bound framework to obtain exact binary solutions.

Table~\ref{tab:random_results} reports the results. 
DCRA achieves nearly a $100\%$ win rate against all three baselines across all tested problem sizes. demonstrating consistent superiority in objective value.
The mean relative differences are negative in all cases, with DCRA achieving on average a 
$53\%$ lower objective value compared to the other methods.

In terms of runtime, DCRA is significantly faster than Gurobi, though slower than MPEC-EPM and L2-box ADMM.
While MPEC-EPM is the fastest among all methods, it consistently yields the worst objective values.

\begin{table*}[ht]
\caption{Performance comparison on randomly generated problems}\label{tab:random_results}
\centering
\resizebox{.9\textwidth}{!}{  
\begin{tabular}{ll|rrr|rrrr}
\hline
\multicolumn{1}{c}{\multirow{2}{*}{n}} & \multicolumn{1}{c|}{\multirow{2}{*}{d}} & \multicolumn{3}{c|}{DCRA's Win Rate   Against:} & \multicolumn{4}{c}{Average CPU Time}         \\ \cline{3-9} 
\multicolumn{1}{c}{}                   & \multicolumn{1}{c|}{}                   & Gurobi       & MPEC-EPM      & L2-box ADMM      & DCRA     & Gurobi   & MPEC-EPM & L2-box ADMM \\ \hline
100                                    & 50                                      & 100\%        & 100\%         & 100\%            & 0.0170   & 0.0737   & 0.0031   & 0.0060      \\
100                                    & 100                                     & 100\%        & 100\%         & 100\%            & 0.0145   & 0.0807   & 0.0034   & 0.0063      \\
100                                    & 200                                     & 100\%        & 100\%         & 100\%            & 0.0203   & 0.1116   & 0.0054   & 0.0090      \\
200                                    & 200                                     & 100\%        & 100\%         & 100\%            & 0.0340   & 0.1638   & 0.0075   & 0.0112      \\
300                                    & 300                                     & 100\%        & 100\%         & 100\%            & 0.0551   & 0.3080   & 0.0114   & 0.0150      \\
300                                    & 500                                     & 100\%        & 100\%         & 100\%            & 0.1218   & 0.6111   & 0.0151   & 0.0183      \\
500                                    & 500                                     & 100\%        & 100\%         & 100\%            & 0.9115   & 2.7283   & 0.1432   & 0.1552      \\
300                                    & 1000                                    & 100\%        & 100\%         & 100\%            & 0.2258   & 1.9576   & 0.0215   & 0.0249      \\
500                                    & 1000                                    & 100\%        & 100\%         & 100\%            & 1.2865   & 8.5689   & 0.2341   & 0.2424      \\
1000                                   & 1000                                    & 100\%        & 100\%         & 100\%            & 1.8787   & 5.4533   & 0.1957   & 0.1989      \\
300                                    & 2000                                    & 100\%        & 100\%         & 100\%            & 0.4629   & 6.8323   & 0.0635   & 0.0720      \\
500                                    & 2000                                    & 100\%        & 100\%         & 100\%            & 1.0622   & 12.8283  & 0.1595   & 0.1609      \\
1000                                   & 2000                                    & 100\%        & 100\%         & 100\%            & 3.1180   & 19.6319  & 0.4112   & 0.4206      \\
2000                                   & 2000                                    & 100\%        & 100\%         & 99\%             & 7.3039   & 29.8193  & 1.4853   & 1.6716      \\
3000                                   & 2000                                    & 100\%        & 100\%         & 88\%             & 13.4699  & 37.2612  & 3.6330   & 3.7831      \\ \hline\hline
\multicolumn{1}{c}{\multirow{2}{*}{n}} & \multicolumn{1}{c|}{\multirow{2}{*}{d}} & \multicolumn{3}{c|}{Mean Relative Difference}   & \multicolumn{4}{c}{Average Objective Value}  \\ \cline{3-9} 
\multicolumn{1}{c}{}                   & \multicolumn{1}{c|}{}                   & Gurobi       & MPEC-EPM      & L2-box ADMM      & DCRA     & Gurobi   & MPEC-EPM & L2-box ADMM \\ \hline
100                                    & 50                                      & -64\%        & -77\%         & -59\%            & $1.44^{02}$  & $3.96^{02}$  & $6.33^{02}$  & $3.53^{02}$   \\
100                                    & 100                                     & -59\%        & -73\%         & -56\%            & $3.30^{02}$  & $8.06^{02}$  & $1.21^{03}$  & $7.54^{02}$    \\
100                                    & 200                                     & -47\%        & -62\%         & -60\%            & $8.50^{02}$  & $1.60^{03}$  & $2.27^{03}$  & $2.12^{03}$    \\
200                                    & 200                                     & -61\%        & -74\%         & -51\%            & $8.89^{02}$  & $2.27^{03}$  & $3.46^{03}$  & $1.81^{03}$    \\
300                                    & 300                                     & -61\%        & -75\%         & -45\%            & $1.62^{03}$  & $4.20^{03}$  & $6.37^{03}$  & $2.97^{03}$    \\
300                                    & 500                                     & -51\%        & -66\%         & -51\%            & $3.43^{03}$  & $6.94^{03}$  & $1.00^{04}$  & $7.02^{03}$    \\
500                                    & 500                                     & -61\%        & -74\%         & -36\%            & $3.52^{03}$  & $8.98^{03}$  & $1.38^{04}$  & $5.54^{03}$    \\
300                                    & 1000                                    & -35\%        & -52\%         & -51\%            & $9.00^{03}$  & $1.38^{04}$  & $1.87^{04}$  & $1.85^{04}$    \\
500                                    & 1000                                    & -46\%        & -62\%         & -48\%            & $9.60^{03}$  & $1.79^{04}$  & $2.54^{04}$  & $1.85^{04}$    \\
1000                                   & 1000                                    & -60\%        & -74\%         & -22\%            & $1.01^{04}$  & $2.52^{04}$  & $3.88^{04}$  & $1.29^{04}$    \\
300                                    & 2000                                    & -22\%        & -38\%         & -38\%            & $2.16^{04}$  & $2.77^{04}$  & $3.47^{04}$  & $3.48^{04}$    \\
500                                    & 2000                                    & -31\%        & -48\%         & -47\%            & $2.46^{04}$  & $3.57^{04}$  & $4.73^{04}$  & $4.67^{04}$    \\
1000                                   & 2000                                    & -46\%        & -62\%         & -37\%            & $2.72^{04}$  & $5.04^{04}$  & $7.20^{04}$  & $4.33^{04}$    \\
2000                                   & 2000                                    & -60\%        & -74\%         & -11\%            & $2.87^{04}$  & $7.15^{04}$  & $1.10^{05}$  & $3.22^{04}$    \\
3000                                   & 2000                                    & -67\%        & -79\%         & -4\%             & $2.92^{04}$  & $8.75^{04}$  & $1.40^{05}$  & $3.06^{04}$     \\ \hline
\end{tabular} 
}

{\small Note: The superscript $k$ denotes a scale of $10^k$. }
\end{table*}

\subsection{Performance on binary compressed sensing}

The binary compressed sensing (BCS) dataset is generated by simulating a linear measurement model
\[
y = A x_0,
\]
where \(x_0 \in \{0,1\}^N\) is a sparse binary vector and \(A \in \mathbb{R}^{M\times N}\) is a random sensing matrix.
Following \cite{doi2024phase}, for a given signal dimension \(N\), we choose the compression ratio \(\alpha = M/N\), the sparsity rate \(\rho\), and the bias parameter \(\mu\) that controls the mean of each matrix entry. 
Each element of \(A\) is independently drawn from a Gaussian distribution
\(
A_{ij} \sim \mathcal{N}\left(\frac{\mu}{N}, \frac{1}{N}\right),
\) 
so that the columns have expected norm of order one. The true signal \(x_0\) is generated by a Bernoulli process with
\(\Pr\{x_{0i}=1\}=\rho\), ensuring approximately \(\rho N\) active components. 

By varying \(\alpha\), \(\rho\), and \(\mu\), one can construct a family of instances spanning different compression and sparsity levels, as done in .
Each instance is fully characterized by (\(N,\alpha,\rho,\mu\)) and the random seed. The original formulation of binary compressed sensing is given by
\begin{equation}\label{eq:ori01}
  \underset{x \in \{0,1\}^N}{\min} \; \|Ax - y\|_1 + \lambda \|x\|_1,
\end{equation}
where \(A \in \mathbb{R}^{M \times N}\) is the sensing matrix, \(y \in \mathbb{R}^M\) is the measurement vector, and \(\lambda > 0\) is a regularization parameter. The term \(\|x\|_1\) promotes sparsity since \(x \in \{0,1\}^N\).

Because our algorithms are designed for binary variables in \(\{-1,1\}\), we apply the variable change \(x = (z + e)/2\), where \(z \in \{-1,1\}^N\). Substituting into~\eqref{eq:ori01} yields an equivalent problem:
\begin{equation}\label{eq:pm1}
  \underset{z \in \{-1,1\}^N}{\min} \; \left\| \tfrac{1}{2} A z - \Big(y - \tfrac{1}{2} A e\Big) \right\|_1
  + \tfrac{\lambda}{2} e^\top z.
\end{equation}
Defining \(\tilde A = \tfrac{1}{2}A\), \(\tilde b = y - \tfrac{1}{2}A e\), and \(c = (\lambda/2)e\), the compact form is
\[
  \underset{z \in \{-1,1\}^N}{\min} \; \|\tilde A z - \tilde b\|_1 + c^\top z.
\]
 
The linear term $c_z^\top z$ cannot in general be expressed purely through $u=A'z-b'$. To fit the template $\min_z \widetilde f(\bar A z - \bar b)$, we introduce an augmented formulation:
\begin{align}\label{eq:aug}
  \bar A &= \begin{bmatrix} A' \\ I_N \end{bmatrix}, &
  \bar b &= \begin{bmatrix} b' \\ 0 \end{bmatrix}, &
  \bar u &= \bar A z - \bar b = \begin{bmatrix} A' z - b' \\ z \end{bmatrix}.
\end{align}
We then define
\begin{equation}\label{eq:f-aug}
  \bar f(\bar u) = \|\bar u_{1:M}\|_1 + c_z^\top \bar u_{M+1:M+N}.
\end{equation}
With this construction,
\[
  \bar f(\bar A z - \bar b) = \|A' z - b'\|_1 + c_z^\top z.
\]
 
The augmented loss $\bar f$ is separable across its two blocks. On the first $M$ coordinates, $\bar f$ reduces to the $\ell_1$-norm. Its proximal operator is soft-thresholding:
  \[
    \operatorname{prox}_{\gamma \|\cdot\|_1}(x) = \mathrm{sign}(x)\cdot \max\{|x|-\gamma,0\}.
  \]
  The corresponding Moreau envelope is
  \[
    \operatorname{env}_{\gamma \|\cdot\|_1}(x) = 
    \sum_{i=1}^M \begin{cases}
      \tfrac{1}{2\gamma} x_i^2, & |x_i|\le \gamma,\\[1ex]
      |x_i|-\tfrac{\gamma}{2}, & |x_i|>\gamma.
    \end{cases}
  \]
On the last $N$ coordinates, $\bar f$ is linear: $c_z^\top u$. Its proximal operator is a simple shift:
  \[
    \operatorname{prox}_{\gamma \langle c_z,\cdot\rangle}(y) = y - \gamma c_z,
  \]
  and the Moreau envelope is
  \[
    \operatorname{env}_{\gamma \langle c_z,\cdot\rangle}(y) 
    = \langle c_z,y\rangle - \tfrac{\gamma}{2}\|c_z\|_2^2.
  \] 
Hence the prox and envelope of $\bar f$ can be computed efficiently by combining these two blocks.

To solve~\eqref{eq:pm1}, we compare algorithms: MPEC-EPM , DCRA, MPEC-ADM, and GUROBI. Each solver produces a vector \(z\), which we map back to \(\{0,1\}^N\) via \(x = (z + e)/2\). For evaluation, we compute the original objective value~\eqref{eq:ori01}, the mean squared error (MSE) \(\|x - x_0\|_2^2/N\), the Hamming distance between \(x\) and \(x_0\), and wall-clock time.

To apply Gurobi to solve \eqref{eq:ori01}, we also need to reformulate it as an equivalent MILP:
{\small
$$
\begin{array}{rc}
\min_{x, z} & e\tp z +\lambda e\tp x\\
{\rm s.t.} & z \ge -(Ax - b), \\
           & z \ge Ax - b, \\
           & x \in \{0, 1\}^N, z\in \R^M.
\end{array}
$$}

\begin{figure*}[ht] 
\centering
\begin{minipage}{.85\linewidth} 
\hspace{-1cm}
\includegraphics[width=1.1\textwidth]{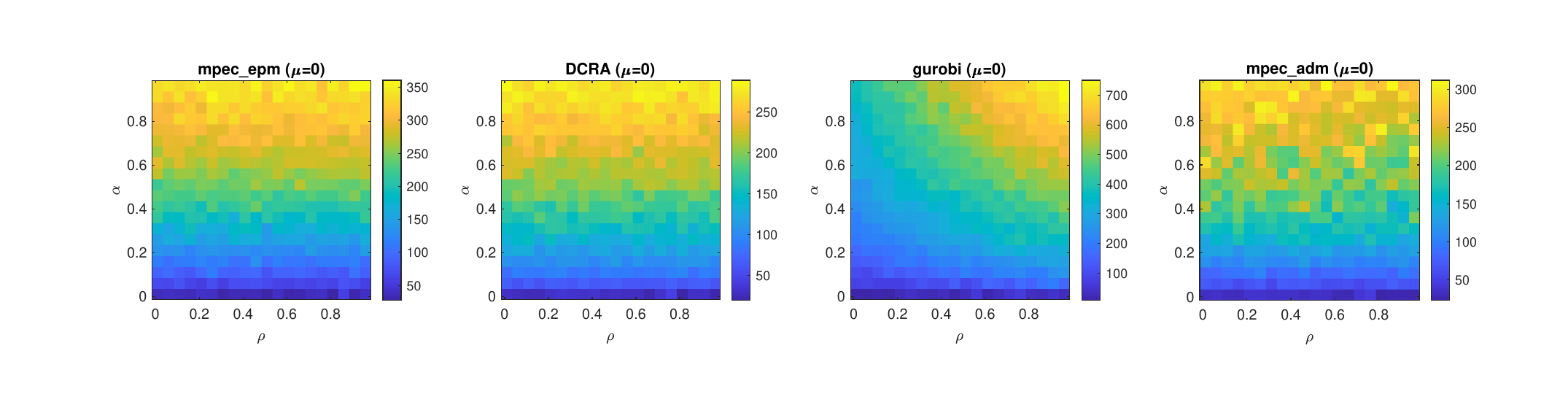}
\vspace{-1.2cm}
\end{minipage}
\begin{minipage}{.85\linewidth}
\hspace{-1cm}
\includegraphics[width=1.1\textwidth]{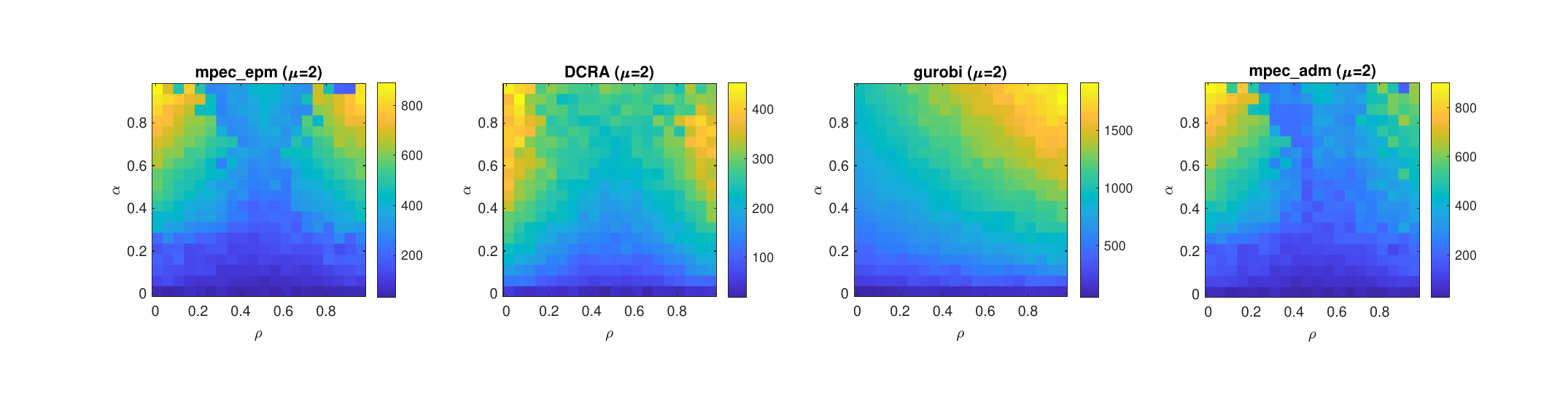}
\vspace{-1.2cm}
\end{minipage}
\begin{minipage}{.85\linewidth}
\hspace{-1cm}
\includegraphics[width=1.1\textwidth]{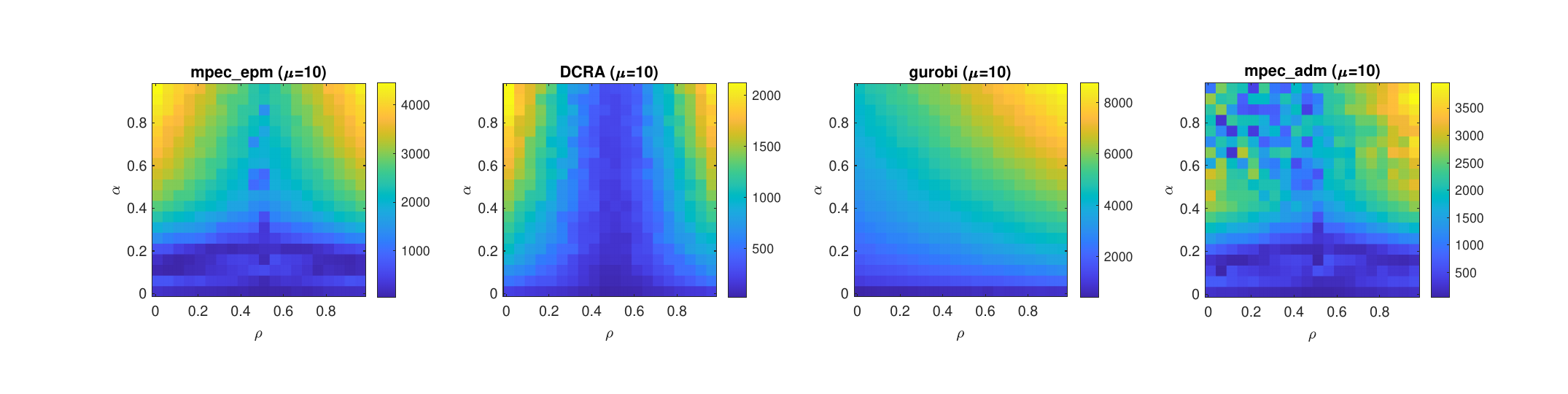} 
\end{minipage}
\caption{\small Phase diagrams for $\mu = 0, 2, 10$ cases. }\label{fig:bcs}
\end{figure*}

Figure \ref {fig:bcs} presents phase diagrams  of objective values versus parameters \( \rho \) and \( \alpha \) for four solvers under three bias levels (\( \mu = 0, 2, 10 \)). Table \ref{tab:bcs} reports average runtime and objective values of four solvers on binary compressed sensing under different biases $\mu$ and sparsity $\rho$. DCRA achieves the lowest objective values consistently, confirming its stronger optimization capability, while MPEC-EPM and MPEC-ADM run fastest.

\hspace{1.5cm}

\vspace{-1.2cm}

\begin{table*}[h]
\centering
\caption{Performance Comparison on Binary Compressed Sensing}\label{tab:bcs}
\resizebox{.95\textwidth}{!}{  
\begin{tabular}{c|l|rrrrrrrrrrr}
\hline
\multicolumn{1}{l|}{$\mu$ = 0}  & $\rho$       & 1\%    & 11\%   & 21\%   & 31\%   & 41\%   & 51\%   & 61\%   & 71\%   & 81\%   & 91\%   \\ \hline
\multirow{4}{*}{Avg. Time}    & mpec\_epm & 0.0739 & 0.0723 & 0.0688 & 0.0706 & 0.0688 & 0.0698 & 0.0719 & 0.0717 & 0.0725 & 0.0710 \\
                             & DCRA      & 0.4965 & 0.4979 & 0.5001 & 0.5012 & 0.5006 & 0.4836 & 0.4981 & 0.5066 & 0.4950 & 0.5084 \\
                             & Gurobi    & 0.4618 & 0.4531 & 0.4459 & 0.4485 & 0.4437 & 0.4599 & 0.4484 & 0.4568 & 0.4535 & 0.4525 \\
                             & mpec\_adm & 0.0751 & 0.0748 & 0.0692 & 0.0721 & 0.0723 & 0.0730 & 0.0734 & 0.0736 & 0.0722 & 0.0724 \\ \hline
\multirow{4}{*}{Avg. Obj. Val.} & mpec\_epm & 217.20 & 219.85 & 215.11 & 216.59 & 221.54 & 213.83 & 216.29 & 218.81 & 216.38 & 221.65 \\
                             & DCRA      & 182.33 & 182.95 & 184.21 & 179.41 & 182.23 & 184.25 & 182.88 & 184.01 & 185.24 & 180.05 \\
                             & Gurobi    & 241.34 & 275.05 & 303.96 & 329.82 & 359.00 & 386.72 & 413.01 & 440.86 & 453.59 & 474.44 \\
                             & mpec\_adm & 197.98 & 189.91 & 192.16 & 190.88 & 192.02 & 189.28 & 185.33 & 190.88 & 185.62 & 179.29 \\ \hline
 
\hline
\multicolumn{1}{l|}{$\mu$ = 2}  & $\rho$       & 0.01   & 0.11   & 0.21   & 0.31   & 0.41   & 0.51   & 0.61     & 0.71     & 0.81     & 0.91     \\ \hline
\multirow{4}{*}{Avg. Time}    & mpec\_epm & 0.0658 & 0.0655 & 0.0654 & 0.0656 & 0.0649 & 0.0651 & 0.0662   & 0.0653   & 0.0640   & 0.0654   \\
                             & DCRA      & 0.4836 & 0.5208 & 0.5019 & 0.5264 & 0.4960 & 0.4901 & 0.5020   & 0.5116   & 0.5004   & 0.4969   \\
                             & Gurobi    & 0.3955 & 0.4035 & 0.3953 & 0.3907 & 0.3939 & 0.3935 & 0.3906   & 0.3855   & 0.3921   & 0.3826   \\
                             & mpec\_adm & 0.0658 & 0.0661 & 0.0666 & 0.0651 & 0.0664 & 0.0666 & 0.0660   & 0.0666   & 0.0659   & 0.0672   \\ \hline
\multirow{4}{*}{Avg. Obj. Val.} & mpec\_epm & 500.00 & 433.93 & 375.50 & 269.88 & 237.53 & 252.43 & 244.96   & 305.22   & 384.34   & 420.96   \\
                             & DCRA      & 315.32 & 271.81 & 227.10 & 201.53 & 186.98 & 182.20 & 190.28   & 203.84   & 229.08   & 280.08   \\
                             & Gurobi    & 641.19 & 711.16 & 759.10 & 833.88 & 892.97 & 975.99 & 1,035.97 & 1,091.03 & 1,158.40 & 1,222.42 \\
                             & mpec\_adm & 503.83 & 438.47 & 378.77 & 261.84 & 204.37 & 233.10 & 227.19   & 252.34   & 322.03   & 344.80   \\ \hline 
                             
\hline
\multicolumn{1}{l|}{$\mu$ = 10} & $\rho$       & 0.01     & 0.11     & 0.21     & 0.31     & 0.41     & 0.51     & 0.61     & 0.71     & 0.81     & 0.91     \\ \hline
\multirow{4}{*}{Avg. Time}    & mpec\_epm & 0.0650   & 0.0652   & 0.0662   & 0.0672   & 0.0651   & 0.0659   & 0.0661   & 0.0652   & 0.0664   & 0.0665   \\
                             & DCRA      & 0.6348   & 0.6440   & 0.6601   & 0.6736   & 0.6848   & 0.7283   & 0.6903   & 0.6648   & 0.6580   & 0.6374   \\
                             & Gurobi    & 0.0801   & 0.0796   & 0.0786   & 0.0823   & 0.0806   & 0.0886   & 0.0903   & 0.0898   & 0.0892   & 0.0906   \\
                             & mpec\_adm & 0.0658   & 0.0660   & 0.0646   & 0.0664   & 0.0658   & 0.0673   & 0.0645   & 0.0664   & 0.0666   & 0.0666   \\ \hline
\multirow{4}{*}{Avg. Obj. Val.} & mpec\_epm & 2,614.78 & 2,319.62 & 2,061.92 & 1,795.14 & 1,574.88 & 1,158.58 & 1,567.82 & 1,825.19 & 2,093.73 & 2,368.07 \\
                             & DCRA      & 1,432.58 & 1,133.98 & 848.24   & 541.57   & 298.18   & 186.64   & 306.21   & 595.97   & 895.64   & 1,190.97 \\
                             & Gurobi    & 2,984.03 & 3,276.34 & 3,569.60 & 3,881.01 & 4,133.50 & 4,471.29 & 4,749.66 & 5,057.28 & 5,365.27 & 5,652.05 \\
                             & mpec\_adm & 1,761.00 & 1,423.93 & 1,337.58 & 1,106.91 & 1,107.50 & 975.05   & 1,440.93 & 1,703.56 & 1,901.69 & 2,177.99 \\ \hline
\end{tabular} 
}
 
\end{table*}

\subsection{Performance on supervised hashing Code Learning}

In this experiment, we evaluate the performance of our method on a supervised discrete hashing problem, a key task in hashing-based retrieval systems. The supervised hashing model follows the formulation in~\cite{Xiong21, gui2017fast, gui2016supervised} with a mean absolute error (MAE) loss:
\begin{equation}\label{eq:sh}
\underset{X \in \{-1,1\}^{n \times r},\; W \in \mathbb{R}^{d \times r}}{\min} \left\| B - W X^\top \right\|_1 + \frac{\delta}{2} \left\| W \right\|_F^2,
\end{equation}
where \( B \in \mathbb{R}^{d \times n} \) is the data matrix, \( W \in \mathbb{R}^{d \times r} \) is a real-valued projection matrix, and \( X \in \{-1,1\}^{n \times r} \) is the binary hash code matrix. This objective is nonconvex due to the binary constraint on \(X\) and the use of the nonsmooth \( \ell_1 \)-norm.

We conduct experiments on three benchmark datasets: MNIST\footnote{MNIST dataset: \url{https://pjreddie.com/projects/mnist-in-csv}}, CIFAR-10\footnote{CIFAR-10 dataset: \url{https://www.cs.toronto.edu/~kriz/cifar.html}}, and MIRFLICKR-25K\footnote{MIRFLICKR-25K dataset: \url{https://press.liacs.nl/mirflickr/mirdownload.html}}.

To solve problem~\eqref{eq:sh}, we adopt an alternating minimization scheme. Specifically, 
for the \(W\)-subproblem, we fix the binary matrix \(X\) and minimize a smoothed version of the objective using gradient descent:
  \begin{equation}\label{eq:sh_sub}
  \underset{W \in \mathbb{R}^{d \times r}}{\min} h_{\mu}\left( B - W (X^{(k)})^\top \right) + \frac{\delta}{2} \|W\|_F^2,
  \end{equation}
  where \( h_\mu(\cdot) \) denotes the Huber loss function with parameter \(\mu > 0\).

Letting \( R = B - W X^\top \), the gradient of the smoothed objective is given by 
\(
\nabla_W \tilde{f}(W) = -H_\mu'(R) X + \delta W,
\)
where \(H_\mu'\) applies the Huber derivative elementwise.

We use a fixed stepsize \( \alpha = 1/L \), where the Lipschitz constant is estimated as 
\(
L = \frac{nr}{\mu} + \delta. 
\)
For the \(X\)-subproblem, we fix \(W\) and solve:
\[
\underset{X \in \{-1,1\}^{n \times r}}{\min} \left\| B - W X^\top \right\|_1,
\]
using algorithms including DCRA, MPEC-EDM, \( L_2 \)-box ADMM, and the subgradient descent method. The parameters of MPEC-EDM and \( L_2 \)-box ADMM are set as default. 

\begin{algorithm}[h!]
\caption{Alternating Minimization for Smoothed Supervised Hashing}
\begin{algorithmic}[1]
\State \textbf{Input:} Data matrix \(B \in \mathbb{R}^{d \times n}\), regularization parameter \(\delta > 0\), initial matrices \(W^{(0)} \in \mathbb{R}^{d \times r}\), \(X^{(0)} \in \{-1,1\}^{n \times r}\)
\For{$k = 0$ to $K-1$}
    \State Fix \(X^{(k)}\), update \(W^{(k+1)}\) by solving~\eqref{eq:sh_sub} via gradient descent with step size \(\alpha = \frac{1}{\frac{nr}{\mu} + \delta}\).
    \State Fix \(W^{(k+1)}\), update \(X^{(k+1)}\) by solving:
    \[
    \underset{X \in \{-1,1\}^{n \times r}}{\min} \left\| B - W^{(k+1)} X^\top \right\|_1.
    \]
\EndFor
\State \textbf{Output:} Final hash code \(X^{(K)}\) and projection \(W^{(K)}\)
\end{algorithmic}
\end{algorithm}


\begin{figure*}[h] 
\centering
\begin{minipage}{.85\linewidth}
\hspace{-2cm}
\includegraphics[width=1.1\textwidth]{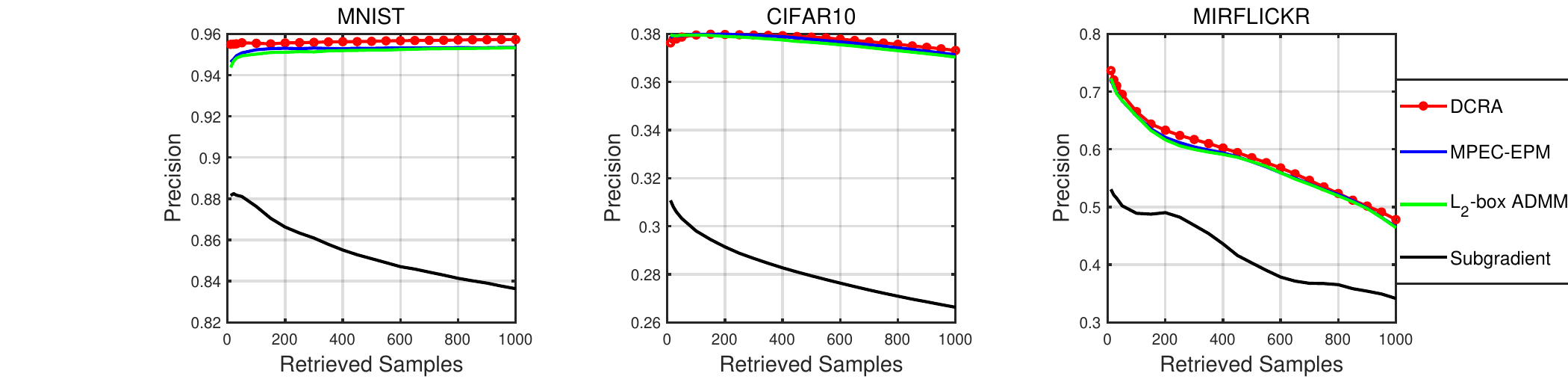}
\vspace{-.3cm}
\caption{\small Precision vs Retrieval Samples}\label{fig:pr}
\end{minipage}
\begin{minipage}{.85\linewidth}
\hspace{-2cm}
\includegraphics[width=1.1\textwidth]{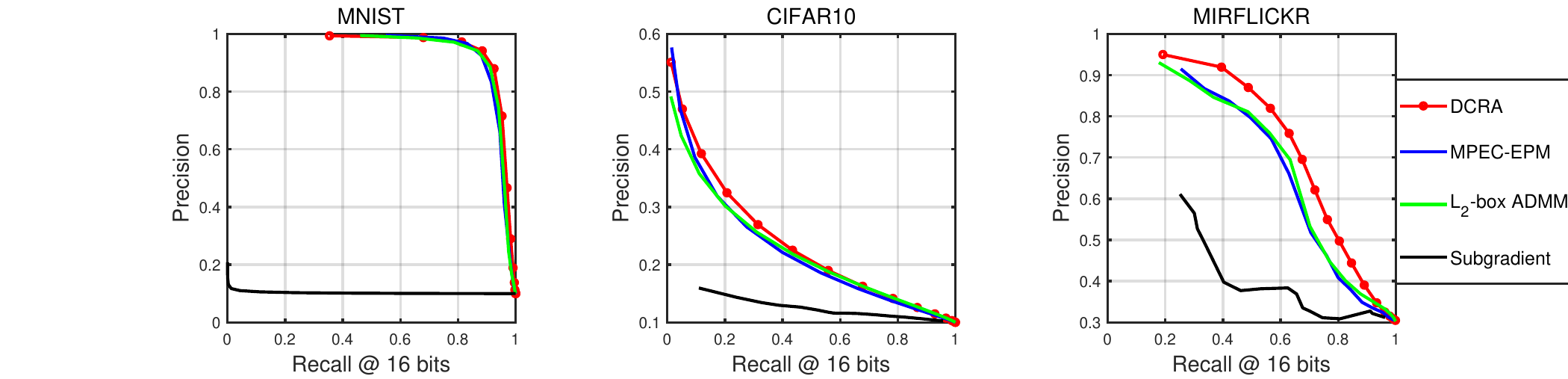}
\vspace{-.3cm}
\caption{\small Precision vs. Recall \@ 16}\label{fig:pr1}
\end{minipage}
\begin{minipage}{.85\linewidth}
\hspace{-2cm}
\includegraphics[width=1.1\textwidth]{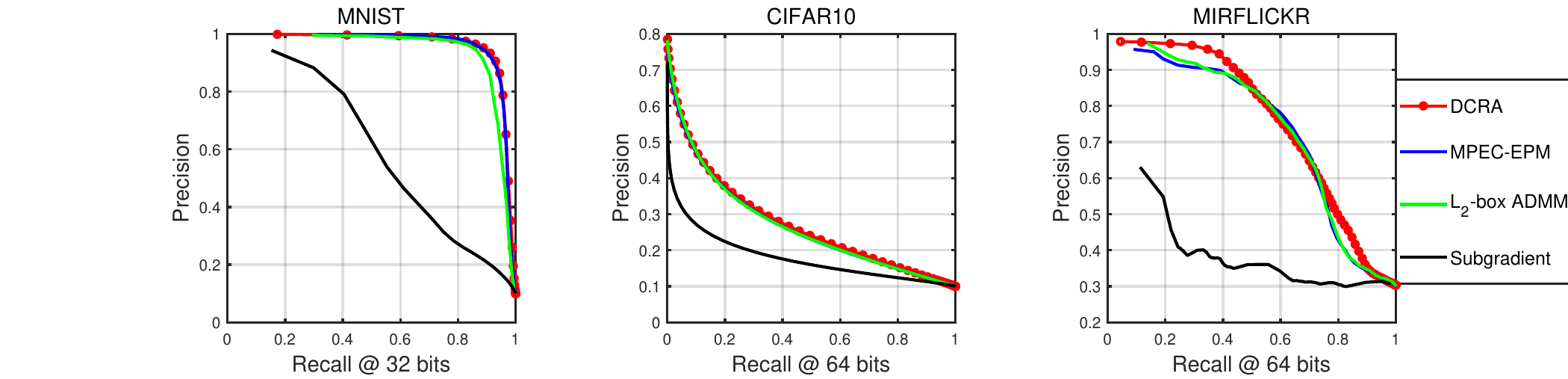}
\vspace{-.3cm}
\caption{\small Precision vs. Recall \@ 32}\label{fig:pr2}
\end{minipage}
\end{figure*}
\vspace{-.5cm}

\begin{figure*}[h] 
\centering
\begin{minipage}{.85\linewidth}
\hspace{-2cm}
\includegraphics[width=1.1\textwidth]{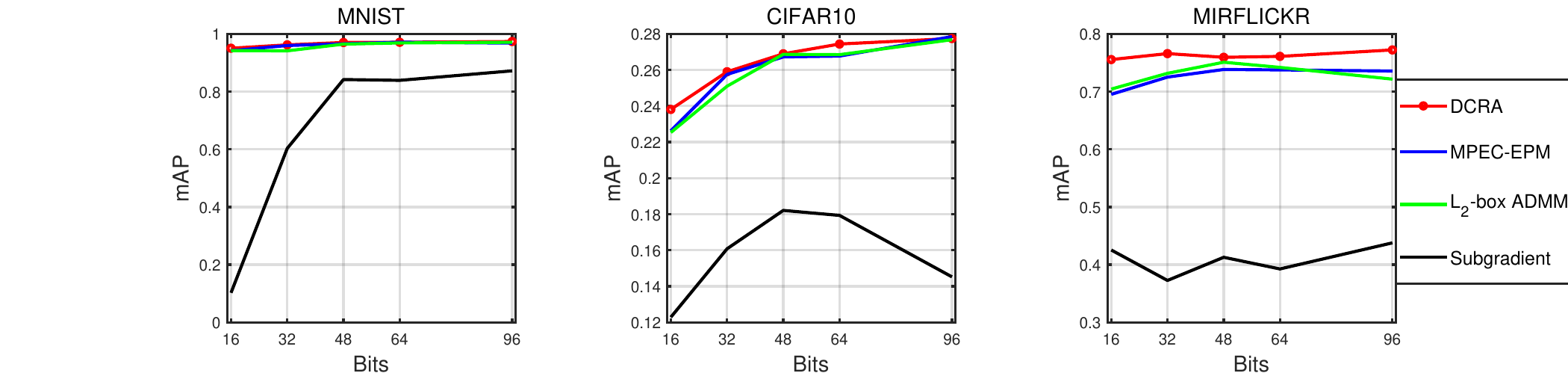}
\vspace{-.3cm}
\caption{\small mAP Curve}\label{fig:map}
\end{minipage}
\begin{minipage}{.85\linewidth}
\hspace{-2cm}
\includegraphics[width=1.1\textwidth]{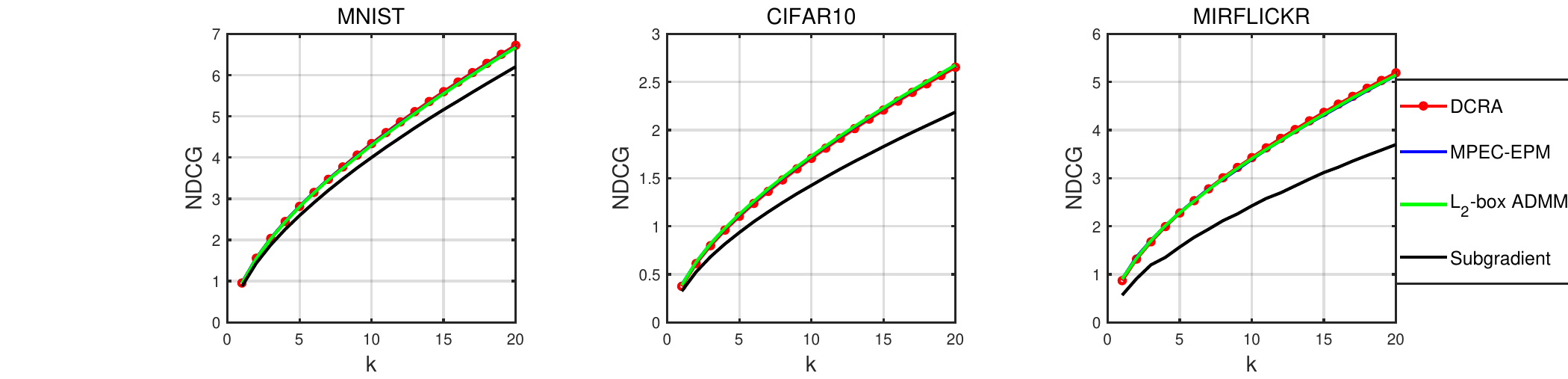}
\vspace{-.3cm}
\caption{\small NDCG with $k = 20$}\label{fig:NDCG}
\end{minipage}
\begin{minipage}{.85\linewidth}
\hspace{-2cm} 
\includegraphics[width=1.1\textwidth]{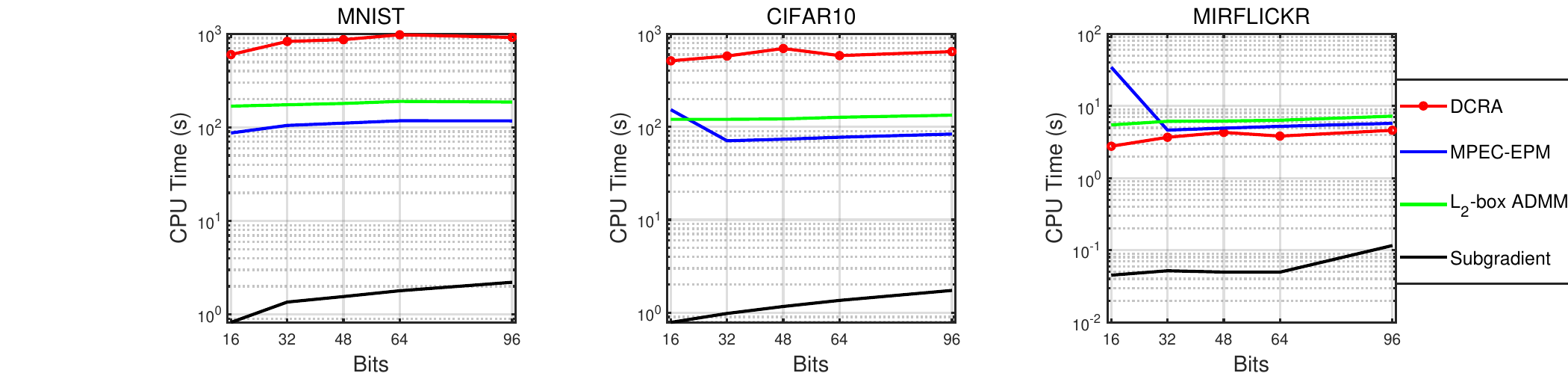}
\vspace{-.3cm}
\caption{\small CPU Time}\label{fig:time}
\end{minipage}
\end{figure*} 

\vspace{.5cm} 
We follow the implementation setup from \cite{gui2017fast, shen2015supervised, Xiong21} and use a linear hash function to encode queries into binary codes. The performance of each method is evaluated using three standard metrics: Mean Average Precision (MAP), Precision–Recall (PR) curve, and Normalized Discounted Cumulative Gain (NDCG).
MAP measures the average precision of retrieved items for each query and provides a single-value summary of ranking quality. The PR curve visualizes the trade-off between precision and recall across different retrieval thresholds. NDCG evaluates the ranking quality by measuring the cumulative gain of retrieved items, discounted by their rank positions, and normalized by the ideal ranking gain. These metrics provide a comprehensive evaluation of retrieval performance. 

Although DCRA incurs higher computational cost compared to the other methods (see Figure~\ref{fig:time}), it consistently outperforms them in all other quality metrics. As shown in Figure~\ref{fig:pr}, DCRA maintains competitive precision across varying retrieval depths. The PR curves for the 16-bit and 32-bit configurations are shown in Figures~\ref{fig:pr1} and \ref{fig:pr2}, respectively. Figure~\ref{fig:map} further demonstrates that DCRA achieves superior MAP scores across all tested datasets in the 16- to 96-bit range. In terms of the NDCG metric (Figure~\ref{fig:NDCG}), DCRA also achieves comparable or better performance across all datasets.

\section{Conclusion}\label{sec:conc}
\IEEEPARstart{I}{n} this work, we developed a novel matrix relaxation approach for solving nonsmooth binary optimization problems, with a focus on absolute value regression. Starting from a rank-constrained SDP reformulation, we introduced a DC surrogate for the rank constraint and derived an exact penalty formulation that enables continuous optimization over a smooth manifold. To solve the resulting penalized problem efficiently, we proposed a matrix factorization-based algorithm (DCRA) with provable convergence guarantees. Our theoretical analysis established the correspondence between the stationary points of the penalized and original formulations, as well as the convergence of the inner iterations under a mild spectral gap condition. Extensive numerical results on synthetic and supervised hashing problems confirmed the effectiveness of our method in producing high-quality binary solutions with favorable computational efficiency. Future work may explore extensions to constrained binary programs, other classes of nonsmooth objectives, and adaptive strategies for tuning the penalty parameters.

\vspace{-0.2cm}
\section*{Acknowledgments}
 This work is supported by the National Natural Science Foundation of China under project No.12371299. 

{\appendices 
\vspace{-0.2cm}
\section{Proof of Proposition \ref{prop:inner-descent}}

\begin{IEEEproof}
\emph{(i)}
By $L_{\tilde f}$-smoothness,
{\small
\[
  \tilde f(V^{k,l+1})\le \tilde f(V^{k,l})
  +\langle \nabla\tilde f(V^{k,l}),\Delta^{k,l}\rangle
  +\tfrac{L_{\tilde f}}{2}\|\Delta^{k,l}\|_F^2,
\]}
where $ \Delta^{k,l}:=V^{k,l+1}-V^{k,l}$. 
By optimality of $V^{k,l+1}$ in \eqref{eq:lin-sub} and feasibility of $V^{k,l}$,
{\small
\begin{align*}
  &\langle \nabla\tilde f(V^{k,l}),\Delta^{k,l}\rangle
  +\rho_k\langle \Gamma^{k,l},\Delta^{k,l}\rangle \\
  &+\rho_k\big(\|V^{k,l+1}\|_F^2-\|V^{k,l}\|_F^2\big)
  +\tfrac{L}{2}\|\Delta^{k,l}\|_F^2
  \ \le\ 0.
\end{align*} }
Since $\|V\|_F^2$ is constant on $\mathcal S$, the term $(\|V^{k,l+1}\|_F^2-\|V^{k,l}\|_F^2)$ cancels. Adding the two above and using $\tilde\psi(V^{k,l+1})\le \tilde\psi(V^{k,l})+\langle \Gamma^{k,l},\Delta^{k,l}\rangle$ gives
{\small
\[
  \Phi_k(V^{k,l+1})-\Phi_k(V^{k,l})
  \ \le\ -\tfrac{L-L_{\tilde f}}{2}\|\Delta^{k,l}\|_F^2.
\]}

\emph{(ii)}
Since $\Phi_k(V^{k,l})$ is nonincreasing and bounded below, it converges. Boundedness of $\{V^{k,l}\}$ holds because $\mathcal S$ is closed and bounded.

\emph{(iii)}
The optimality condition of \eqref{eq:lin-sub} gives, for some $W^{k,l+1}\in N_{\mathcal S}(V^{k,l+1})$,
\[
  0 =\nabla\tilde f(V^{k,l})+\rho_k\Gamma^{k,l}+2\rho_k V^{k,l+1}
  +L\Delta^{k,l}+W^{k,l+1}.
\]
Hence,
\begin{align*}
  &\nabla\tilde f(V^{k,l+1})+2\rho_k V^{k,l+1}+\rho_k\Gamma^{k,l}+W^{k,l+1}\\
  = &\big(\nabla\tilde f(V^{k,l+1})-\nabla\tilde f(V^{k,l})\big)-L\Delta^{k,l}.
\end{align*} 
Taking norms and using $L_{\tilde f}$-Lipschitz continuity of $\nabla\tilde f$ yields \eqref{eq:prop4}. 
\end{IEEEproof}

\vspace{-0.2cm}
\section{Proof of Theorem \ref{thm:inner-complexity} }

\begin{IEEEproof} 
Let $\Delta^{k,l}:=V^{k,l+1}-V^{k,l}$. By Proposition~\ref{prop:inner-descent}(iii),  
taking minimum on both sides of \eqref{eq:prop4} over $l=0,\dots,T-1$ gives
\begin{equation}\label{eq:min-direct}
\underset{0\le l<T}{\min}\mathcal G\big(V^{k,l+1}, \Gamma^{k,l}\big)
\ \le\ C \underset{0\le l<T}{\min} \|\Delta^{k,l}\|_F.
\end{equation}
For any nonnegative sequence $\{a_l\}$ we have
$\min_{0\le l<T} a_l \le \sqrt{\tfrac{1}{T}\sum_{l=0}^{T-1} a_l^2}$. Applying this to
$a_l=\|\Delta^{k,l}\|_F$ yields
\begin{equation}\label{eq:min-le-rms}
\underset{0\le l<T}{\min}\|\Delta^{k,l}\|_F
\ \le\ \sqrt{\frac{1}{T}\sum_{l=0}^{T-1}\|\Delta^{k,l}\|_F^2}.
\end{equation}
By Proposition~\ref{prop:inner-descent}(i),
\begin{equation}\label{eq:sum-steps}
\sum_{l=0}^{T-1}\|\Delta^{k,l}\|_F^2
  \le\ \frac{2}{L-L_{\tilde f}}\Big(\Phi_k(V^{k,0})-\Phi_k(V^{k,T})\Big) .
\end{equation}
 Since \eqref{eq:lin-sub} is a majorization of \eqref{eq:smooth-fact}, we have for all $V\in\mathcal S$, 
\(
  \Phi_{\delta}^\star
  \;=\;\inf_{V\in\mathcal S}\Phi_\delta(V)
  \;\le\; \Phi_\delta(V)
  \;\le\;\Phi_k(V),
\)
and hence
\(
  \Phi_k(V^{k,0})-\Phi_k(V^{k,T})
  \;\le\;
  \Phi_k(V^{k,0})-\Phi_{\delta}^\star.
\) 
Combining \eqref{eq:min-direct}, \eqref{eq:min-le-rms}, and \eqref{eq:sum-steps} gives \eqref{eq:itr_cmplx},  
where $V^{k} = V^{k,0}$ by Step \ref{alg:step1} of Algorithm \ref{AlgA}. 
In particular, if
\(
T\ \ge\ \left\lceil \frac{2C^2}{L-L_{\tilde f}}\,
\frac{\Phi_k(V^{k})-\Phi_{\delta}^\star}{\epsilon^2}\right\rceil,
\)
then $\min_{0\le l<T}\mathcal G(\cdot)\le \epsilon$, i.e., the inner loop attains an
$\epsilon$-stationary point in $\mathcal O(\epsilon^{-2})$ iterations.

For any $\epsilon_v>0$, choosing
$T\ge \frac{\Phi_k(V^{k})-\Phi_{\delta}^\star}{c\,\epsilon_v^{\,2}}$ guarantees that there is an index $\ell<T$ with
$\|V^{k,\ell+1}-V^{k,\ell}\|_F\le \epsilon_v$, so Algorithm \ref{AlgA} terminates. 
\end{IEEEproof}

\vspace{-0.2cm}

\section{Proof of Lemma \ref{lem:complexity}}

\begin{IEEEproof}
Write $V^{k,0}:=V^{l}$ and $U^{k,0}:=V^{k,0}\in\mathcal S$.
Take
\(
\Gamma^{k,0} := - 2\,V^{k,0}\,P_1 P_1^\top  \in \partial \widetilde\psi(V^{k,0}),
\)
where $P_1\in\mathbb R^p$ is the leading eigenvector of $(V^{k,0})^\top V^{k,0}$.
For each $j\in[p]$, the $j$-th column of $\Gamma^{k,0}$ is
\begin{align*}
&\Gamma^{k,0}_j \;=\; -\,2\,V^{k,0}P_1\, (P_1)_j \Longrightarrow \\
&
\|\Gamma^{k,0}_j\|  =  2 \|V^{k,0}P_1\| |P_{j1}|
 =  2\,\|V^{k,0}\|\,|P_{j1}|
 \ge  2 |P_{j1}|,
\end{align*}
where we used that $P_1$ is a right singular vector associated with $\|V^{k,0}\|$ and
$\|V^{k,0}\|\ge \|V^{k,0}_j\|=1$ for all $j$. Hence
\begin{equation}\label{eq:Gamma-col-lb}
\|\Gamma^{k,0}_j\|\ \ge\ 2\mu\; \text{for all }j\in[p],
\;\mu:=\min_{1\le j\le p}|P_{j1}|>0.
\end{equation}
Following \eqref{eq:lin-sol}, define
\begin{equation}\label{eq:G-def}
G^{k,0} :=\frac{1}{L+2\rho_k} \left( L V^{k,0} - \nabla \widetilde{f}(V^{k,0}) - \rho_k \Gamma^{k,0} \right).
\end{equation}
Using the Lipschitz continuity
of $\nabla\widetilde f$ with constant $L_{\widetilde f}$ and
$\|V^{k,0}\|_F=\sqrt p$, $\|E/\sqrt m\|_F=\sqrt p$, we obtain
{\small
\begin{align}\label{eq:varpi}
\big\|LV^{k,0}&-\nabla\widetilde f(V^{k,0})\big\|_F
 \le L\|V^{k,0}\|_F \notag \\
&+ L_{\widetilde f}\|V^{k,0}-\tfrac1{\sqrt m}E\|_F  + \big\|\nabla\widetilde f(\tfrac1{\sqrt m}E)\big\|_F
 \le \varpi,
\end{align}
}
with $\varpi:=(L+2L_{\widetilde f})\sqrt p+\|\nabla\widetilde f(\tfrac1{\sqrt m}E)\|_F$.
Taking $j$-th columns in \eqref{eq:G-def}, the triangle inequality, and
\eqref{eq:Gamma-col-lb}–\eqref{eq:varpi} yield
{\small
\begin{align}\label{eq:col-norm-lb}
\|G^{k,0}_j\|& \ge \frac{1}{L+2\rho_k}\lf(\rho_k\|\Gamma^{k,0}_j\|
-\big\|[LV^{k,0}-\nabla\widetilde f(V^{k,0})]_j\big\|\rg)\notag \\
& \ge \frac{2\rho_k\mu-\varpi}{L+2\rho_k}.
\end{align}}
Under $\rho_k\ge\overline\rho=(Lp+\varpi)/(\mu c_0)$ with $c_0\in(0,1)$ we have
$\rho_k\mu\ge \varpi$ and therefore from \eqref{eq:col-norm-lb}
\begin{equation}\label{eq:col-norm-lb2}
\|G^{k,0}_j\|\ \ge\ \frac{\rho_k\mu}{L+2\rho_k}\ >\ 0\qquad\forall j.
\end{equation} 
Let $D:=\mathrm{Diag} \big(1/\|G^{k,0}_1\|,\dots,1/\|G^{k,0}_p\|\big)$.
Then $V^{k,1}=\operatorname{Proj}_{\mathcal S}(G^{k,0})=G^{k,0}D$. From \eqref{eq:col-norm-lb2}
\begin{equation}\label{eq:D-bounds}
\|D\|\ \le\ \frac{L+2\rho_k}{\rho_k\mu}, 
\|D\|_F \le \sqrt p\,\|D\| \le \sqrt p\,\frac{L+2\rho_k}{\rho_k\mu}.
\end{equation}
Define $\overline{\Gamma}^{k,0}:= -\tfrac12\,\Gamma^{k,0}D$, which has rank one according to the definition of $\Gamma^{k,0}$.
Using \eqref{eq:G-def} and $V^{k,1}=G^{k,0}D$,
\(
V^{k,1}-\overline{\Gamma}^{k,0}=\frac{1}{L+2\rho_k}\lf(\frac{L}{2 }\Gamma^{k,0}D
\ +(LV^{k,0}-\nabla\widetilde f(V^{k,0})) D\rg).
\)
Taking Frobenius norms and using
$\|\Gamma^{k,0}\|\le 2\|V^{k,0}\|\le 2\|V^{k,0}\|_F=2\sqrt p$ together with
\eqref{eq:varpi} and \eqref{eq:D-bounds}, we obtain
{\footnotesize
\begin{align*}
\operatorname{dist}&(V^{k,1},\overline{\mathcal R})
 \le \|V^{k,1}-\overline{\Gamma}^{k,0}\|_F\\
& \le \frac{1}{L+2\rho_k}\lf(\frac{L}{2 }\,\|\Gamma^{k,0}\|\,\|D\|_F
 + \|LV^{k,0}-\nabla\widetilde f(V^{k,0})\|_F\,\|D\|\rg) \\[-1mm]
&\ \le\  \frac{1}{L+2\rho_k}\lf(L\sqrt{p} \cdot
\frac{L+2\rho_k}{\rho_k\mu}\ + \varpi \cdot
\frac{L+2\rho_k}{\rho_k\mu}\rg)
\ =\ \frac{Lp+\varpi}{\rho_k\mu}.
\end{align*}}
With $\rho_k\ge\overline\rho=(Lp+\varpi)/(\mu c_0)$ this yields
$\operatorname{dist}(V^{k,1},\overline{\mathcal R})\le c_0$.
Finally, by the Eckart–Young–Mirsky theorem, define $\overline{\mathcal R}:=\{V\in\mathbb{R}^{m\times p}\mid 
\mathrm{rank}(V)\le 1\}$, 
\begin{align*}
\operatorname{dist}(V^{k,1},\overline{\mathcal R})^2
& =\ \|V^{k,1}\|_F^2 - \sigma_1(V^{k,1})^2\\
& =\ \|V^{k,1}\|_F^2 - \|V^{k,1}\|^2.
\end{align*}
\end{IEEEproof}

\vspace{-0.3cm}
\section{Proof of Proposition \ref{prop:inner-progress}}

\begin{IEEEproof}
(i) Fix $l\in\mathbb N$ with $\frac{\varepsilon}{2}\le p-\|V^{k,l}\|^2\le c_0^2$.
By optimality of $V^{k,l+1}$ for the linearized subproblem,
$\widetilde\Phi_k (V^{k,l+1};V^{k,l})\le \widetilde\Phi_k (V;V^{k,l})$ for any $V\in\mathcal S$, hence
{\small
\begin{align*}
\rho_k\!&\left\langle \Gamma^{k,l}, V^{k,l+1} - V\right\rangle 
\le \left\langle \nabla \widetilde f (V^{k,l}),\, V - V^{k,l+1}\right\rangle \\[-1mm]
& \qquad + \frac{L}{2}\!\left(\|V - V^{k,l}\|_F^2 - \|V^{k,l+1} - V^{k,l}\|_F^2\right)  \\[-1mm]
&\le 2\sqrt p\!\left(L_{\widetilde f}\,\big\|V^{k,l}-\tfrac{1}{\sqrt m}E\big\|_F
+\big\|\nabla\widetilde f\!\left(\tfrac1{\sqrt m}E\right)\big\|_F\right) \\[-1mm]
&\qquad + \frac{L}{2}\|V\|_F^2 + L\,\|V-V^{k,l+1}\|_F\,\|V^{k,l}\|_F \\[-1mm]
&\le 2\sqrt p\!\left(L_{\widetilde f}\cdot 2\sqrt p
+\big\|\nabla\widetilde f\!\left(\tfrac1{\sqrt m}E\right)\big\|_F\right)
+ \frac{L}{2}p + 2Lp \\[-1mm]
&= 4L_{\widetilde f}p + 2\sqrt p\,\big\|\nabla\widetilde f\!\left(\tfrac1{\sqrt m}E\right)\big\|_F
+ \tfrac{5}{2}Lp :=\,\varpi_{\mathrm{lin}},
\end{align*}}
where we used Lipschitz continuity of $\nabla\widetilde f$ and the facts $\|V\|_F=\|V^{k,l}\|_F=\sqrt p$ and $\|V-V^{k,l+1}\|_F\le 2\sqrt p$ for $V\in\mathcal S$.
Rearranging gives, for any $V\in\mathcal S$,
{\small
\begin{align}\label{eq:spectral}
- \left\langle \Gamma^{k,l}, V\right\rangle 
&\le \rho_k^{-1}\varpi_{\mathrm{lin}} + \left\langle\Gamma^{k,l}, V^{k,l+1}\right\rangle \notag\\
&\le \rho_k^{-1}\varpi_{\mathrm{lin}}+ \|\Gamma^{k,l}\|_*\;\|V^{k,l+1}\| \notag\\
& = \rho_k^{-1}\varpi_{\mathrm{lin}}+ 2\|V^{k,l}\|\;\|V^{k,l+1}\|,
\end{align}}
where we used our spectral choice $\Gamma^{k,l}=-2V^{k,l}u_1u_1^\top$ and the rank-one identity $\|\Gamma^{k,l}\|_*=2\|V^{k,l}\|$.
Now take the SVD $V^{k,l}=Q\,\mathrm{Diag}(\sigma)\,U^\top$, with top singular vectors $q_1,u_1$.
Let $\widehat u_1\in\mathbb R^p$ be defined by $(\widehat u_1)_j=u_{j1}/|u_{j1}|$ and set $\widehat V:=q_1\widehat u_1^\top\in\mathcal S$.
Then 
\(
-\langle \Gamma^{k,l},\,\widehat V\rangle
=2\,\langle V^{k,l}u_1u_1^\top,\,q_1\widehat u_1^\top\rangle
=2\,\|V^{k,l}\|\sum_{j=1}^p |u_{j1}|.
\) 
Insert $V=\widehat V$ into \eqref{eq:spectral} to obtain
\begin{equation}\label{eq:key-prod}
2\,\|V^{k,l}\|\sum_{j=1}^p |u_{j1}|
\ \le\ \rho_k^{-1}\varpi_{\mathrm{lin}}+2\,\|V^{k,l}\|\,\|V^{k,l+1}\|.
\end{equation}
From $\|V^{k,l}_j\|=1$ for each column and the SVD identity,
\begin{align*}
\sigma_1(&V^{k,l})^2\,u_{j1}^2
=1-\sum_{i\ge2}\sigma_i(V^{k,l})^2\,u_{ji}^2 \\[-1mm]
&\ge 1-\sum_{i\ge2}\sigma_i(V^{k,l})^2
= 1-\big(\|V^{k,l}\|_F^2-\|V^{k,l}\|^2\big).
\end{align*}
Since $p-\|V^{k,l}\|^2=\|V^{k,l}\|_F^2-\|V^{k,l}\|^2\le c_0^2$, we have
\begin{equation}\label{eq:S1}
|u_{j1}|\ \ge\ \frac{\sqrt{1-c_0^2}}{\|V^{k,l}\|}\quad\forall j.
\end{equation} 
By $\sum_{j=1}^p u_{j1}^2=1$, there exists an index $\hat{\jmath}$ with $u_{\hat{\jmath}1}^2\le 1/p$.
Using $x\ge x^2$ for $x\in[0,1]$,
{\small
\begin{equation}\label{eq:bound}
\sum_{j=1}^p |u_{j1}|
\ \ge\ |u_{\hat{\jmath}1}|+\sum_{j\ne\hat{\jmath}}u_{j1}^2
\ \ge\ \frac{\sqrt{1-c_0^2}}{\|V^{k,l}\|}+\Big(1-\frac{1}{p}\Big).
\end{equation}}
Plugging \eqref{eq:bound} into \eqref{eq:key-prod} yields
{\small
\[
2\,\|V^{k,l}\|\,\|V^{k,l+1}\|
\ \ge\ 2\sqrt{1-c_0^2}\ +\ 2\Big(1-\frac{1}{p}\Big)\|V^{k,l}\|\ -\ \rho_k^{-1}\varpi_{\mathrm{lin}}.
\]}
Since $p-\|V^{k,l}\|^2\le c_0^2$ implies $\|V^{k,l}\|\ge \sqrt{p-c_0^2}$, using $\,\|V^{k,l+1}\|^2 +\|V^{k,l}\|^2\,\ge 2\|V^{k,l}\|\,\|V^{k,l+1}\|$, it gives 
{\small
\[ 
\|V^{k,l+1}\|^2\ \ge\ \|V^{k,l}\|^2\ +\  2\Big(\,\Big(1-\tfrac{1}{p}\Big)\sqrt{p-c_0^2}\ +\ \sqrt{1-c_0^2}\,\Big) - \frac{\varpi_{\mathrm{lin}}}{\rho_k}\,. 
\]}
Let $\eta(c_0,p):= 2\Big(1-\tfrac{1}{p}\Big)\sqrt{p-c_0^2}\ +\ 2\sqrt{1-c_0^2}$. Thus, choosing $\rho_k\ge \dfrac{2\,\varpi_{\mathrm{lin}}}{\eta(c_0,p)}$ guarantees \eqref{eq:gap_v}.

(ii)
Let $g_\ell:=p-\|V^{k,\ell}\|^2$. 
Assume for contradiction that $g_\ell>\varepsilon$ for all $\ell=1,\dots,\widehat l:=\big\lceil 2(c_0^2-\varepsilon)/\eta(c_0,p)\big\rceil+1$.
Then for $\ell=1,\dots,\widehat l-1$ we have $g_\ell\in[\varepsilon/2,\,c_0^2]$, so item~(i) applies to get
\(
\|V^{k,\ell+1}\|^2 \ge \|V^{k,\ell}\|^2+\tfrac12\,\eta(c_0,p)
 \Longleftrightarrow 
g_{\ell+1} \le g_\ell-\tfrac12\,\eta(c_0,p).
\)
Summing these $\widehat l - 1$ inequalities yields $g_{\widehat l}\le c_0^2-\tfrac12(\widehat l-1)\eta(c_0,p)\le\varepsilon$, a contradiction.
Hence there exists $1\le \overline l\le \widehat l$ with $g_{\overline l}\le\varepsilon$.

By Lemma~\ref{lem:complexity} and $\rho_k\ge\overline\rho$, we have $g_1\le c_0^2$.
If $g_1\le\varepsilon$, take $\overline l=1$ and we are done. Otherwise $\varepsilon<g_1\le c_0^2$.
 
If $g_j\in[\varepsilon/2,\varepsilon]$, applying item~(i) again gives
$g_{j+1}\le g_j-\tfrac12\,\eta(c_0,p)\le \varepsilon$.
If $g_j<\varepsilon/2$, use \eqref{eq:spectral} with $V=V^{k,j}$ to obtain
\begin{align*}
2\,\|V^{k,j}\| & \|V^{k,j+1}\|\ \ge\ 2\,\|V^{k,j}\|^2 - \rho_k^{-1}\,\varpi_{\mathrm{lin}}\\
\Longrightarrow &\|V^{k,j+1}\|^2  \ge \|V^{k,j}\|^2 - \rho_k^{-1}\,\varpi_{\mathrm{lin}}.
\end{align*}
Thus $g_{j+1}\le g_j+\rho_k^{-1}\varpi_{\mathrm{lin}}
< \varepsilon/2+\varepsilon/2=\varepsilon$ since $\rho_k\ge 2\,\varpi_{\mathrm{lin}}/\varepsilon$.
Therefore, $g_l\le\varepsilon$ for all $l\ge\overline l$.

(iii)
By \eqref{eq:S1} with $l=1$ and repeating the same argument at every $l\ge1$, we have 
{\small
\begin{equation}\label{eq:S2}
|u_{1j}|
 \ge \frac{\sqrt{1-c_0^2}}{\sigma_1}
 =  \frac{\sqrt{1-c_0^2}}{\|V^{k,l}\|}
 \ge  \frac{\sqrt{1-c_0^2}}{\sqrt{p}}
 > 0,
 \;\forall j\in[p].
\end{equation}}
By the fact that $P_1(V^{k,l}) = u_1(V^{k,l})$, the leading eigenvector $P_1(V^{k,l})$ of $(V^{k,l})^\top V^{k,l}$ has no zero entries for every $l\ge1$.
\end{IEEEproof}

\section{Proof of Theorem \ref{thm:feas-obj}}

\begin{IEEEproof} 
Define $\Delta^{k,l}:= V^{k,l+1} - V^{k,l}$. For a fixed $k$, the update \eqref{eq:lin-sol} with $L\ge L_{\widetilde f}$ yields
{\small
\begin{equation}\label{eq:bound-norm}
\rho_k  \|V^{k,l+1}\|_F^2 
\le 
\rho_k \|V^{k,l}\|_F^2 
- \langle \nabla \widetilde{f}(V^{k,l})+\notag 
 \rho_k\Gamma^{k,l}, \Delta^{k,l} \rangle  - \frac{L}{2} \|\Delta^{k,l}\|_F^2.
\end{equation}}
Combining with the descent lemma for  $\widetilde f$,
\(
\widetilde f(V^{k,l+1})-\rho_k\|V^{k,l+1}\|^2
\ \le\ \widetilde f(V^{k,l})-\rho_k\|V^{k,l}\|^2. 
\)
Thus, 
$\widetilde f(\cdot)-\rho_k\|\cdot\|^2$ is nonincreasing. Since $V^{k+1}\in \{V^{k,l}\}$ and $V^{k,0}=V^k$, this gives
\begin{equation}\label{eq:inner-descent-k}
\widetilde f(V^{k+1})-\rho_k\|V^{k+1}\|^2\ \le\ \widetilde f(V^{k})-\rho_k\|V^{k}\|^2.
\end{equation}
Adding $(\rho_k-\rho_{k+1})\|V^{k+1}\|^2$ to both sides of \eqref{eq:inner-descent-k} and telescoping from $k=k^*$ to $\overline{k}-1$ yields
\begin{align}\label{eq:outer-telescope-k}
\widetilde f(V^{\overline{k}})-\rho_{\overline{k}}\|V^{\overline{k}}\|^2
&\le \widetilde f(V^{k^*})-\rho_{k^*}\|V^{k^*}\|^2 \notag \\[-2mm]
& + \sum_{j=k^*}^{\overline{k}-1}\big(\rho_j-\rho_{j+1}\big) \|V^{j+1}\|^2.
\end{align} 
Write the SVD of $V^{\overline{k}} = Q{\rm Diag}(\sigma_1, \dots, \sigma_r)U\tp$ and  eigen-decomposition of $X^{\overline{k}}= U{\rm Diag}(\lambda_1, \dots, \lambda_r, 0,\dots, 0)U\tp$. Since $X^{\overline{k}}=(V^{\overline{k}})^\top V^{\overline{k}}=\sum_{i=1}^p \lambda_i P_iP_i^\top$ and $x^{\overline{k}}:=\|V^{\overline{k}}\|\,P_1$, we have $x^{\overline{k}}(x^{\overline{k}})^\top=\lambda_1\,P_1P_1^\top$. Then, 
{\small
\[
\big\|X^{\overline{k}}-x^{\overline{k}}(x^{\overline{k}})^\top\big\|_F
=\Big(\sum_{i=2}^p \lambda_i^2\Big)^{\frac{1}{2}}\\
\le \sum_{i=2}^p \sigma_i^2
= \|V^{\overline{k}}\|_F^2-\|V^{\overline{k}}\|^2\le\varepsilon. 
\]}
Taking diagonals gives 
$\|x^{\overline{k}}\!\circ x^{\overline{k}}-e\| \le \|x^{\overline{k}}(x^{\overline{k}})^\top-X^{\overline{k}}\|_F \le \varepsilon$.
Since $\hat{f}$ is Lipschitz on $\Omega$ with modulus $L_{\hat f}$ and $X^{\overline{k}}\in\Omega$,
\begin{align}\label{eq:prevs}
\hat{f}(X^{\overline{k}})
&\ge\ \hat{f}\!\big(x^{\overline{k}}(x^{\overline{k}})^\top\big)
 - L_{\hat f} \big\|X^{\overline{k}}-x^{\overline{k}}(x^{\overline{k}})^\top\big\|_F \notag\\[-1mm]
& \ge\ \hat{f}\!\big(x^{\overline{k}}(x^{\overline{k}})^\top\big) - L_{\hat f}\,\varepsilon.
\end{align} 
Since $\|V^{k^*}\|\ge \|V^{k^*}\|_F/\sqrt{\mathrm{rank}(V^{k^*})} = \sqrt{p/r^*}$, using $\widetilde f(V^{\overline{k}})=\hat{f}\!\big((V^{\overline{k}})^\top V^{\overline{k}}\big)=\hat{f}(X^{\overline{k}}))$ and \eqref{eq:outer-telescope-k}, we obtain
{\small
\begin{align*}
\hat{f}(X^{\overline{k}}) - \widetilde f(V^{k^*})
&\le\ \rho_{\overline{k}}\|V^{\overline{k}}\|^2\ -\ \rho_{k^*}\,\frac{p}{r^*}\\[-1mm]
& +\ \sum_{j=k^*}^{\overline{k}-1}\big(\rho_j-\rho_{j+1}\big)\,\|V^{j+1}\|^2.
\end{align*} }
Combine above with \eqref{eq:prevs} to get \eqref{eq:obj-gap-dcra}. By Proposition \ref{prop:inner-progress}(ii), set $g_{\overline{k}}:=p-\|V^{\overline{k}}\|^2\in[0,\epsilon]$, so that
$\|V^{\overline{k}}\|^2=p-g_{\overline{k}}$.
Since $\rho_{j+1}\ge\rho_j$ and $\|V^{j+1}\|^2\ge 1$ for all $j$,
{\small
\[
\sum_{j=k^*}^{\overline{k}-1}\big(\rho_j-\rho_{j+1}\big)\,\|V^{j+1}\|^2
\ \le\ \sum_{j=k^*}^{\overline{k}-1}\big(\rho_j-\rho_{j+1}\big)
\ =\ \rho_{k^*}-\rho_{\overline{k}}.
\]}
Thus,
{\footnotesize
\begin{align*}
\hat{f}\!\big(x^{\overline{k}}(x^{\overline{k}})^\top\big)-\widetilde f(V^{k^*})
&\le \rho_{\overline{k}}\,\|V^{\overline{k}}\|^2
-\rho_{k^*}\frac{p}{r^*}
+\big(\rho_{k^*}-\rho_{\overline{k}}\big)
+L_{\hat f}\,g_{\overline{k}}\\
&= \rho_{\overline{k}}\,(p-g_{\overline{k}})
-\rho_{k^*}\frac{p}{r^*}
+\rho_{k^*}-\rho_{\overline{k}}
+L_{\hat f}\,g_{\overline{k}}\\
&= \rho_{\overline{k}}(p-1)
+\rho_{k^*}\Big(1-\frac{p}{r^*}\Big)
+\big(L_{\hat f}-\rho_{\overline{k}}\big)\,g_{\overline{k}}\\
&\le \rho_{\overline{k}}(p-1)
+\rho_{k^*}\Big(1-\frac{p}{r^*}\Big)
+\big(L_{\hat f}-\rho_{\overline{k}}\big)\,\epsilon.
\end{align*} }
\end{IEEEproof}

}

\section{Biography Section}

\vspace*{-.5cm}
\begin{IEEEbiographynophoto}{Lianghai Xiao}
received the Ph.D. degree from School of Mathematics, University of Birmingham, UK in 2020. He was a postdoctoral researcher with the School of Mathematics, South China University of Technology, China. He is currently a lecturer with College of Information Science and Technology, Jinan University, China.
\end{IEEEbiographynophoto}

\vspace*{-.5cm}
\begin{IEEEbiographynophoto}{Yitian Qian}
received the B.S. and the Ph.D. degree from the School of Mathematics, South China University of Technology, China, in 2018 and 2023, respectively.
He is currently a postdoctoral researcher with Department of Data Science and Artificial Intelligence, The Hong Kong Polytechnic University, Hong Kong, China.
\end{IEEEbiographynophoto}

\vspace*{-.5cm}
\begin{IEEEbiographynophoto}{Shaohua Pan}
	received the Ph.D. degree from School of Mathematics, Dalian University of Technology, China, in 2003. She is currently a professor with the School of Mathematics, South China University of Technology, China. Her main research interest involves low-rank and sparsity optimization, as well as nonconvex and nonsmooth optimization.
\end{IEEEbiographynophoto}

\vfill

\end{document}